\documentclass[10pt]{article} 
\usepackage{amsmath,amssymb,amscd}
\usepackage{array}
\usepackage{euscript}
\usepackage{epsfig}
\usepackage{geometry}

\newtheorem{Lemma}{Lemma}[section]
\newtheorem{Theorem}[Lemma]{Theorem}

\newtheorem{Remark}[Lemma]{Remark}

\def\qed{\hspace*{11.5cm} \framebox[0.22cm][r]}

\def\Hom{\op{Hom}}
\def\End{\op{End}}
\def\Aut{\op{Aut}}

\def\spa{\op{span}}

\def\GL{\op{GL}}

\def\Mat{\op{Mat}}
\def\Diag{\op{Diag}}
\def\diag{\op{diag}}
\def\rk{\op{rank}}
\def\R{\mathbb{R}}

\newcommand\op{\operatorname}

\newcommand{\D}{\displaystyle}

\begin{document}
\newgeometry{left=4.7cm,right=4.7cm, top=3.5cm, bottom=4.5cm}
\begin{center}\Large{\bf Indecomposable symplectic $k(C_2\times C_2)$--modules and their quadratic forms}\end{center}
\begin{center}
Lars Pforte \footnote{Department of Mathematics, National University of Ireland, Maynooth, Co. Kildare, Ireland, lars.pforte@nuim.ie} and John Murray \footnote{Department of Mathematics, National University of Ireland, Maynooth, Co. Kildare, Ireland, john.murray@nuim.ie}
\end{center}
\enlargethispage{0.5cm}
\let\thefootnote\relax\footnotetext{Date: 30 Nov 2017}

\begin{center}{\bf Abstract}\end{center}
\small
For the Klein-Four Group $G$ and a perfect field $k$ of characteristic two we determine all indecomposable symplectic $kG$-modules, that is, $kG$-modules with a symplectic, $G$-invariant form which do not decompose into smaller such modules, and classify them up to isometry. Also we determine all quadratic forms that have one of the above symplectic forms as their associated bilinear form and describe their isometry classes.

\section{Introduction}

Let $k$ be a field and let $G$ be a group. The classification of the $kG$-modules up to isomorphism is equivalent to classifying the conjugacy classes of subgroups isomorphic to $G$ in the general linear groups over $k$. In this paper we are interested in $kG$-modules with a $G$-invariant symplectic form. Classifying such modules is equivalent to the classification of the conjugacy classes of subgroups isomorphic to $G$ in the symplectic groups over $k$. We also study $G$-invariant quadratic forms and classifying them up to isometry is equivalent to classifying the conjugacy classes of subgroups isomorphic to $G$ in the orthogonal groups over $k$. In \cite{Wall} the conjugacy classes in symplectic and orthogonal groups have been studied, thus providing an answer to the above problems in the case of cyclic groups. In this paper we focus on the Klein-Four Group. 

Let $G$ be the Klein-Four Group and $k$ be a perfect field of characteristic two. We classify all indecomposable symplectic $kG$-modules up to isometry, that is, $kG$-modules with a $G$-invariant symplectic form and which do not decompose into an orthogonal sum of smaller symplectic $kG$-modules. In particular if $k$ is finite we give the number of different isometry classes of indecomposable symplectic $kG$-modules. Furthermore for each indecomposable symplectic $kG$-module we determine if there exists a $G$-invariant quadratic form whose associated bilinear form coincides with the given symplectic form. Finally we also study the isometry classes of those quadratic forms. 

In chapter two we give a brief introduction to symmetric, symplectic and quadratic forms and define what we mean by an indecomposable symplectic $kG$-module. The indecomposable modules for the Klein-Four Group are fully described in \cite{Conlon} and they can be enumerated using the notation
\[\text{$k_G$, $kG$, $A_n$, $B_n$, $C_n(f)$ and $C_n(\infty)$,}\]
for all $n\geq 1$ an integer and all $f\in k[T]$ an irreducible polynomial in the indeterminate $T$. We describe these modules in more detail in chapter three and the following $kG$-modules give rise to an indecomposable symplectic $kG$-module:
\[(k_G)^2,\ kG,\ (kG)^2,\ A_n\oplus B_n,\ C_n(f),\ C_n(f)^2,\ C_n(\infty),\ C_n(\infty)^2,\]
for all $n\geq 1$ and all irreducible $f\in k[T]$. Here $M^2$ denotes $M\oplus M$. In chapter three we deal with these modules in this order. In the following we give a brief and rather loose summary of our main results. The full results for each module can be found in the Theorems \ref{Theorem_(k_G)2}, \ref{Theorem_D}, \ref{Theorem_D^2}, \ref{Theorem_AnoplusBn}, \ref{Theorem-C_n(f)}, \ref{Theorem-C_n(f)+C_n(f)*}, \ref{Theorem-C_n(infty)} and \ref{Theorem-C_n(infty)+C_n(infty)*}, respectively. Throughout the following let $\EuScript{P}$ be a full set of representatives for the distinct cosets of the additive subgroup $\{x^2+x:\ x\in k\}$ in $(k,+)$.

The case of $(k_G)^2$ is independent of the group and has been studied before. The indecomposable symplectic $kG$-modules can be enumerated by the set $k^*$ and they form one isometry class. Given such a symplectic form, the set of all corresponding $G$-invariant quadratic forms can be enumerated by $k$, while the isometry classes can be enumerated by $\EuScript{P}$.

The case $kG$ is also well-understood. If, as in this paper $G$ is the Klein Four Group, then for every triplet $(b,c,d)\in k^3$ such that $d\neq b+c$ there is a symplectic module. Those triplets where $b+c+d=1_k$ represent the various isometry classes. Finally given a symplectic form the set of all corresponding $G$-invariant quadratic forms can be enumerated by $k$ and all these quadratic forms are isometric.

The isometry classes of the indecomposable symplectic forms on $(kG)^2$ are represented by the paired module and two disjoint sets of forms each enumerated by $k$. That means in the finite case we have $2|k|+1$ different isometry classes. Furthermore for each indecomposable symplectic form the set of all corresponding $G$-invariant quadratic forms can be enumerated by $k^2$ and all these quadratic forms are isometric.

In the case of $A_n\oplus B_n$ the indecomposable symplectic $kG$-module can be enumerated by $\EuScript{F}:=\{(\omega_1,\ldots,\omega_{2n})\in k^{2n}: \omega_1=\ldots=\omega_{j-1}=0_k, \omega_j=1_k, j=1,\ldots,2n+1\}$. In particular in the finite case there are $|k|^{2n-1}+\ldots+|k|+2$ different isometry classes. Furthermore the class associated with the label $(\omega_1,\ldots,\omega_{2n})\in\EuScript{F}$ has a corresponding $G$-invariant quadratic form if and only if $\omega_i=\omega_{i+1}$, for all even $i\in\{2,\ldots,2n-2\}$. In the case of the paired module there is one isometry class of $G$-invariant quadratic forms. Otherwise the isometric classes can be enumerated by $\EuScript{P}$. 

For $C_n(f)$ let $m$ denote the degree of $f$ and let $l$ be the largest integers less than or equal to $\frac{n}{2}$. Furthermore let $K=k[\epsilon]$ be a field extension of $k$ of degree $m$ where $\epsilon$ is a root of $f$. Then the indecomposable symplectic $kG$-modules can be enumerated by $\{(a_1,\ldots,a_l)\in K^l\}$, that is, in the finite case we have $|k|^{ml}$ isometry classes. Finally depending on $m$, $n$ and $\epsilon$ at most one symplectic form has a corresponding $G$-invariant quadratic form. In this case those quadratic forms can be enumerated by $\EuScript{P}$.

For $(C_n(f))^2$ let $m$ and $K$ be as before. Then the isometry classes of all indecomposable symplectic $kG$-modules can be enumerated by the quadruples $(r,s,\varphi,\psi)$, where $r$ and $s$ are integers with $2\leq r\leq n+1$ and $0\leq s\leq \frac{n-r+1}{2}$, $\varphi\in K^s$ and $\psi\in\{(\psi_1,\ldots,\psi_t)\in K^t: \psi_1\neq 1_k\}$, where $t$ denotes the largest integer less than or equal to $\frac{n-r-2s+1}{2}$. Note that $(n+1,0,\emptyset,\emptyset)$ is a permissible label. The existence of corresponding $G$-invariant quadratic forms depend on $r$, $m$, $\epsilon$ and $\varphi$ but not on $s$ or $\psi$. Moreover in the case of the paired module there is one isometry class of corresponding quadratic forms, while otherwise the isometry classes can again be enumerated by $\EuScript{P}$.

Finally if $g_1$  and $g_2$ denote generators of $G$, then $g_1$ and $g_2$ act on $C_n(\infty)$ as $g_2$ and $g_1$ act on $C_n(f)$, respectively, where $f(T)=T\in k[T]$. Consequently our results on indecomposable symplectic $kG$-modules and corresponding $G$-invariant quadratic forms for $C_n(\infty)$ and $C_n(\infty)^2$ are the same as for $C_n(T)$ and $C_n(T)^2$, respectively.

\section{Remarks on Forms and important notation}\label{Chapter2}
{\bf{Remarks on bilinear forms}:} Let $k$ be a field and let $(M,B)$ be a $k$-vector space with bilinear form $B:M\times M\rightarrow k$. Also we identify $B$ with its Gram matrix with respect to a fixed basis $\EuScript{B}$ on $M$. For $x,y\in M$ we represent the corresponding row-vectors over $k$ with respect to $\EuScript{B}$ by $x$ and $y$, respectively. Then $B(x,y)=x^TBy$. Vice versa every $B\in\Mat_n(k)$ defines a bilinear form $B:(x,y)\mapsto x^TBy$ on $M$, where $n$ is the $k$-dimension of $M$.

The maps $\Phi:x\mapsto B(x,\cdot)$ and $\Psi:y\mapsto B(\cdot,y)$ give homomorphisms from $M$ into the dual $M^*=\Hom_k(M,k)$ of $M$. The sets $\{x\in M:\ B(x,y)=0, \text{for all $y\in M$}\}$ and $\{y\in M:\ B(x,y)=0, \text{for all $x\in M$}\}$ are their respective kernels, called left and right radical of $B$, respectively. The form $B$ is called \emph{reflexive} if for all $x,y\in M$ we have that $B(x,y)=0$ implies $B(y,x)=0$. In this case the left and right radical of $B$ coincide and we just talk about the the \emph{radical} of $B$. Observe that $\Phi$ is a $k$-isomorphism if and only if $\Psi$ is a $k$-isomorphism. In this case the radical of $B$ is zero and $B$ is called \emph{non-degenerate}. Furthermore $B$ is non-degenerate if and only if the matrix $B$ is invertible. 

Next let $G$ be a group. If $M$ is a $kG$-module, we say $B$ is $G$-invariant, if $B(gx,gy)=B(x,y)$, for all $x,y\in M$ and all $g\in G$. Note that if $B$ is non-degenerate and $G$-invariant, then $\Phi$ and $\Psi$ are $kG$-isomorphisms. In particular a $kG$-module with a non-degenerate, $G$-invariant bilinear form is self-dual. On the other hand every self-dual $kG$-module $M$ has a non-degenerate, $G$-invariant bilinear form $B$, where $B(x,y):=\varphi(x)(y)$, for all $x,y\in M$, and $\varphi: M\rightarrow M^*$ is a $kG$-isomorphism.

Set $\EuScript{E}:=\End_{kG}(M)$. For every $\alpha\in \EuScript{E}$, let $\alpha$ also denote the corresponding matrix over $k$ with respect to the basis $\EuScript{B}$. Two $G$-invariant forms $B$ and $B'$ are called \emph{isometric} if there is some invertible $\alpha\in \EuScript{E}$ such $B(\alpha(x),\alpha(y))=B'(x,y)$, for all $x,y\in M$. This is equivalent to $\alpha^TB\alpha=B'$. Clearly isometry is an equivalence relation.

Let $N$ be a $kG$-submodule of $M$ and let $B$ be a reflexive bilinear form on $M$. Then $N^{\perp}:=\{y\in M:\ B(x,y)=0\text{, for all $x\in N$}\}$ is also a $kG$-submodule of $M$. Two $kG$-submodules $N_1,N_2$ of $M$ are called \emph{orthogonal} if $B(x,y)=0$, for all $x\in N_1$ and $y\in N_2$, i.e. $N_2\subseteq N_1^{\perp}$ and hence also $N_1\subseteq N_2^{\perp}$. If $M$ is a direct sum of pairwise orthogonal submodules $N_1,\ldots,N_t$, then $M$ is an \emph{orthogonal sum} of $N_1,\ldots,N_t$, and we write $(M,B)=(N_1,B)\perp\ldots\perp (N_t,B)$, or more briefly, $M=N_1\perp\ldots\perp N_t$. If $B$ is non-degenerate on both $M$ and $N$, then $M=N\perp N^{\perp}$. Finally let $N_1,N_2$ be $kG$-modules with respective bilinear forms $B_1$ and $B_2$. Then $(N_1,B_1)\perp(N_2,B_2)$ denotes the module $N_1\oplus N_2$ with the form $(B_1\perp B_2)(x_1+x_2,y_1+y_2):=B_1(x_1,y_1)+B_2(x_2,y_2)$, for all $x_1,y_1\in N_1$ and $x_2,y_2\in N_2$.

A form $B$ is called \emph{symmetric} if $B(x,y)=B(y,x)$, for all $x,y\in M$, and \emph{skew-symmetric} if $B(x,y)=-B(y,x)$, for all $x,y\in M$. Note that a symmetric form is reflexive and corresponds to a symmetric matrix. Vice versa a symmetric matrix gives rise to a symmetric form. A form $B$ is called \emph{alternating} if $B(x,x)=0$, for all $x\in M$. In this case the corresponding matrix $B$ is \emph{hollow}, that is, $B$ has a zero main diagonal. Generally a hollow matrix does not correspond to an alternating form as seen in the example of the paired module discussed below. Furthermore an alternating form $B$ is skew-symmetric, as $B(x+y,x+y)=0$ implies $B(x,y)=-B(y,x)$, for all $x,y\in M$. A form $B$ that is alternating and non-degenerate is called \emph{symplectic}. Finally we call $(M,B)$ a \emph{symmetric/symplectic} $kG$-module if $B$ is symmetric/alternating, $G$-invariant and non-degenerate. In particular every symmetric/symplectic $kG$-module is self-dual.

Next we say a symmetric/symplectic module $(M,B)$ is \emph{indecomposable} if $M$ is not an orthogonal sum of two smaller symmetric/symplectic $kG$-modules. By \cite{GowWillems} we know that for an indecomposable symmetric/symplectic $kG$-module $(M,B)$, either $M$ is indecomposable as a $kG$-module or $M\cong N\oplus N^*$, for some indecomposable $kG$-module $N$. Finally an indecomposable symplectic module is even-dimensional.

Now let $B$ be $G$-invariant and non-degenerate. For $\alpha\in \EuScript{E}$ we define $B_{\alpha}(x,y):=B(\alpha(x),y)$, for all $x,y\in M$. Then $B_{\alpha}$ is a $G$-invariant form on $M$ and $B_{\alpha}$ corresponds to $\alpha^TB$. Thus for every $G$-invariant bilinear form $B'$ on $M$, we have $B'=B_{\alpha}$, where $\alpha$ corresponds to $(B'B^{-1})^T$. Therefore $\{B_{\alpha}:\ \alpha\in \EuScript{E}\}$ gives the space of all $G$-invariant bilinear forms on $M$. 

We define the transpose of $B$ as $B^T(x,y):=B(y,x)$, for $x,y\in M$. Then $B^T=B_{\sigma}$, for some $\sigma\in \EuScript{E}$. Then $\sigma=B^{-T}B$. Next recall that for every $\alpha\in \End_k(M)$ there is some $\alpha^{\circ}\in \End_k(M)$, called \emph{adjoint} of $\alpha$, such that $B(x,\alpha(y))=B(\alpha^{\circ}(x),y)$, for all $x,y\in M$. In fact $\alpha^{\circ}=(B\alpha B^{-1})^T$. Also if $\alpha\in \EuScript{E}$, then $\alpha^{\circ}\in \EuScript{E}$. Observe that now $\sigma^{\circ}=(B(B^{-T}B)B^{-1})^T=B^{-1}B^T=\sigma^{-1}$.

\enlargethispage{0.5cm}
Next consider the $kG$-module $M^*\oplus M=\{(f,m):f\in M^*, m\in M\}$. For $f\in M^*$ and $m\in M$ we set $q(f,m):=f(m)$. Then $q$ is a quadratic form on $M^*\oplus M$. For its associated bilinear form $P:=B_q$ we get $P((f_1,x_1),(f_2,x_2))=f_1(x_2)+f_2(x_1)$, 
for all $f_1,f_2\in M^*$ and $x_1,x_2\in M$. Observe that $P$ is $G$-invariant and non-degenerate. Finally $P$ is symmetric in general and symplectic in characteristic two. We call the module $(M^*\oplus M, P)$ the \emph{paired module}.

\begin{Lemma}\label{Lemma-PairedModules}
Let $M_1,M_2$ be $kG$-modules and let $(M_1^*\oplus M_1,P_1)$, $(M_2^*\oplus M_2,P_2)$ and $((M_1\oplus M_2)^*\oplus (M_1\oplus M_2),P)$ be various paired modules. Then $P$ is isometric to $P_1\perp P_2$.
\end{Lemma}

\begin{Lemma}\label{LemmaMB_perp_M2P_iso_MB3}
Let $k$ have characteristic two and let $(M,B)$ be a symmetric $kG$-module. Then $(M,B)\perp (M^2,P)$ is isometric to $(M,B)\perp (M,B)\perp (M,B)$.
\end{Lemma}

Lemmas \ref{Lemma-PairedModules} and \ref{LemmaMB_perp_M2P_iso_MB3} follow quickly by looking at the respective Gram matrices. Next let $M$ be self-dual with $G$-invariant, non-degenerate bilinear form $B$. Then $P_B((x_1,x_2),(y_1,y_2)):=B(x_1,y_2)+B(y_1,x_2)$, for $x_1,x_2,y_1,y_2\in M$, defines a form on $M^2=M\oplus M$. Note that $(M^2,P_B)$ is a symmetric $kG$-module, which is symplectic in characteristic two. If furthermore $\{b_1,\ldots,b_n\}$ is a basis for $M$, then with respect to the basis $\{(b_1,0),\ldots,(b_n,0),(0,b_1),\ldots,(0,b_n)\}$ for $M^2$ the Gram matrix of $P_B$ equals $\begin{pmatrix} 0&B\\B^T&0\end{pmatrix}$. Finally recall that $M\cong_{kG}M^*$ via the map $x\mapsto B(x,\cdot)$. As $P_B((x_1,x_2),(y_1,y_2))=P((B(x_1,\cdot),x_2),(B(y_1,\cdot),y_2))$ we see that $P_B$ is the paired module. 

Next recall that for an indecomposable $kG$-module $M$ the endomorphism ring $\EuScript{E}=\End_{kG}(M)$ is local, that is, every element in $\EuScript{E}$ is either invertible or nilpotent. The Jacobson radical $J(\EuScript{E})$ is the ideal of all nilpotent elements in $\EuScript{E}$ and $\EuScript{D}:=\EuScript{E}/J(\EuScript{E})$ is a division ring. Since $M\rightarrow M:x\mapsto \lambda\cdot x$ lies in $\EuScript{E}$, for all $\lambda\in k$, it follows that $k$ is isomorphically embedded in $\EuScript{D}$. Note that if $k$ is finite, then $\EuScript{D}$ is a field and if $k$ is algebraically closed then $\EuScript{D}\cong k$.

\begin{Lemma}\label{LemmaM2P_indecomposable}
Let $k$ have characteristic two and let $M$ be an indecomposable and self-dual $kG$-module such that $\EuScript{E}/J(\EuScript{E})$ is a field. Then the paired module $(M^2,P)$ is indecomposable, that is, every submodule of $(M^2,P)$ that is isomorphic to $M$ is degenerate. 
\end{Lemma}

Proof: Let $B$ be a non-degenerate $G$-invariant form on $M$. Then $P_B$ is the form on the paired module. Next let $\phi: M\rightarrow M^2$ be an injection. Then there are $\alpha_1,\alpha_2\in \EuScript{E}$ such that $\phi(x)=(\alpha_1(x),\alpha_2(x))\in M^2$, for all $x\in M$. Let $x,y\in M$. Now
\begin{align*}
P(\phi(x),\phi(y))&=P_B(\phi(x),\phi(y))= 
B(\alpha_1(x),\alpha_2(y))+B(\alpha_1(y),\alpha_2(x))\\&=B(\alpha_2^{\circ}\alpha_1(x),y)+B(\alpha_1^{\circ}\sigma\alpha_2(x),y)=B_{\alpha_2^{\circ}\alpha_1+\alpha_1^{\circ}\sigma\alpha_2}(x,y)
\end{align*}
Hence $(\phi(M),P)$ is isometric to $(M,B_{\alpha_2^{\circ}\alpha_1+\alpha_1^{\circ}\sigma\alpha_2})$. Write $\sigma=\mu 1_M+\delta$, for $\mu\in k$ and $\delta\in J(\EuScript{E})$. As $\delta$ is nilpotent we have
\[\mu1_M+\delta^{\circ}=\sigma^{\circ}=\sigma^{-1}=\mu^{-1}\sum_{n=0}^{\infty}(-\mu^{-1}\delta)^n.\]
It follows that $\mu=\mu^{-1}$ and so $\sigma=1_M+\delta$, that is, $\sigma\equiv 1_M\mod (J(\EuScript{E}))$. Also note that $\alpha_i^{\circ}\equiv \alpha_i\mod(J(\EuScript{E}))$, for $i=1,2$, and $\alpha_1\alpha_2\equiv \alpha_2\alpha_1\mod(J(\EuScript{E}))$. Thus $\alpha_2^{\circ}\alpha_1+\alpha_1^{\circ}\sigma\alpha_2\equiv 0\mod(J(\EuScript{E}))$. In particular $B_{\alpha_2^{\circ}\alpha_1+\alpha_1^{\circ}\sigma\alpha_2}$ is degenerate, as claimed.\\
\qed\\\\

Observe that Lemmas \ref{LemmaMB_perp_M2P_iso_MB3} and \ref{LemmaM2P_indecomposable} imply that in characteristic two for an indecomposable symplectic $kG$-module $(M,B)$ there are two very different decompositions of $(M,B)\perp (M,B)\perp (M,B)$ into indecomposable symplectic modules. In particular symplectic modules do not satisfy the Krull-Schmidt Theorem, as has already been noticed for instance in \cite{WillemsThesis}.

As the main focus of this paper is on characteristic two in the following we wish to point out some more observations that are relevant in even characteristic. Earlier we have seen that an alternating form $B$ on $M$ is skew-symmetric. Hence in characteristic two $B$ is symmetric and thus the corresponding matrix is both hollow and symmetric. Vice versa every hollow and symmetric matrix $B\in\Mat_n(k)$ gives rise to an alternating form on $M$. Assume otherwise. Then there must be a change of basis $X\in\GL_n(k)$ such that $X^TBX$ is not hollow. However this is contradicted by the Lemma \ref{Lemma-XTAX_i,i=...} below. In particular with respect to a fixed basis the alternating forms on $M$ are given precisely by the set of all hollow and symmetric $n\times n$--matrices and the symplectic forms on $M$ are given precisely by the set of all hollow, symmetric and invertible $n\times n$--matrices. The proof of the following lemma is straight forward and omitted here. 

\begin{Lemma}\label{Lemma-XTAX_i,i=...}
In characteristic two let $n,m\geq 1$ be integers. Also let $A\in\Mat_{n}(k)$ and $X\in\Mat_{n,m}(k)$ be matrices such that $A$ is symmetric. Then for all $i\in\{1,\ldots,m\}$ 
\[(X^TAX)_{i,i}=\sum_{s=1}^n A_{s,s}\cdot (X_{s,i})^2.\]
Consequently if $A$ is hollow then $X^TAX$ is hollow. In particular if $X$ is invertible, then $A$ is hollow if and only if $X^TAX$ is hollow.
\end{Lemma}

\begin{Remark}\label{Remark-symplectic_forms_on_M^2}
In characteristic two let $M$ be an indecomposable self-dual $kG$-module with basis $\{e_1,\ldots,e_n\}$. Also assume that $\EuScript{E}/J(\EuScript{E})$ is a field. We face this situation in sections \ref{subsection-kG+kG} and \ref{subsection-Cn(f)^2} and so we briefly wish to discuss $G$-invariant symplectic forms on $M^2=M\oplus M^*$. With respect to the basis $\{(e_1,0),\ldots,(e_n,0),(0,e_1),\ldots,(0,e_n)\}$ of $M^2$, let $S=\begin{pmatrix} X&Y\\W&Z\end{pmatrix}\in\GL_{2n}(k)$ be a $G$-invariant symplectic form on $M^2$, where $X,Y,Z,W\in\Mat_n(k)$. Note that the restriction of $S$ to $\spa\{(e_1,0),\ldots,(e_n,0)\}$ is $X$ and the restriction of $S$ to $\spa\{(0,e_1),\ldots,(0,e_n)\}$ is $Z$. Hence both $X$ and $Z$ correspond to $G$-invariant alternating forms on $M$. Also if $X$ or $Z$ is non-degenerate, then $(M^2,S)$ decomposes. Hence if $(M^2,S)$ is indecomposable, then $X$ and $Z$ are degenerate. Furthermore $W=Y^T$,  as $S$ is symmetric. If we define $B(m_1,m_2):=S((m_1,0),(0,m_2))$, for $m_1,m_2\in M$, then $B$ is a $G$-invariant form on $M$ corresponding to $Y$. Hence for every indecomposable symplectic $kG$-module $(M^2,S)$ we have that $S$ corresponds to an invertible matrix $\begin{pmatrix} X&Y\\Y^T&Z\end{pmatrix}$, where $X,Z\in\Mat_n(k)$ correspond to $G$-invariant alternating degenerate forms and $Y\in\Mat_n(k)$ corresponds to a $G$-invariant form. 

On the other hand every such matrix clearly gives rise to a $G$-invariant symplectic form $S$ on $M^2$. Furthermore we claim that for such $S$ the symplectic $kG$-module $(M^2,S)$ is indecomposable. Let $B$ be the form as in the previous paragraph, that is, $Y$ corresponds to $B$. Moreover there are $\alpha,\gamma\in \EuScript{E}:=\End_{kG}(M)$ such that $X$ and $Z$ correspond to $B_{\alpha}$ and $B_{\gamma}$, respectively. Also note that since $X$ and $Z$ are singular, we have $\alpha,\gamma\in J(\EuScript{E})$. Next assume that $(M^2,S)$ is decomposable. Then there is an injection $\phi: M\rightarrow M^2$ such that $(\phi(M),S)$ is a symplectic $kG$-module. Note that there are $\alpha_1,\alpha_2\in\EuScript{E}$ such that $\phi(x)=(\alpha_1(x),\alpha_2(x))$, for all $x\in M$. Also recall that there is $\sigma\in\EuScript{E}$ such that $B^T=B_{\sigma}$. Overall we get
\begin{align*}
S(\phi&(x),\phi(y))=S((\alpha_1(x),\alpha_2(x)),(\alpha_1(y),\alpha_2(y)))\\&=B_{\alpha}(\alpha_1(x),\alpha_1(y))+B(\alpha_1(x),\alpha_2(y))+B_{\sigma}(\alpha_2(x),\alpha_1(y))+B_{\gamma}(\alpha_2(x),\alpha_2(y))\\
&=B_{\tau}(x,y),
\end{align*}
where $\tau=\alpha_1^{\circ}\alpha\alpha_1+\alpha_2^{\circ}\alpha_1+\alpha_1^{\circ}\sigma\alpha_2+\alpha_2^{\circ}\gamma\alpha_2$. Now as in the proof of Lemma \ref{LemmaM2P_indecomposable} we get $\tau\equiv 0\mod{J(\EuScript{E})}$, that is, $B_{\tau}$ is degenerate, which is a contradiction. In particular $(M^2,S)$ is indecomposable. 
\end{Remark}

{\bf{Remarks on quadratic forms}:} 
A map $q:M\rightarrow k$ is called a \emph{quadratic form} on $M$ if $q(a x)=a^2q(x)$, for all $x\in M$ and all $a\in k$, and $B_q(x,y):=q(x+y)-q(x)-q(y)$ is a bilinear form. We call $(M,q)$ a \emph{quadratic space} and $B_q$ the \emph{associated bilinear form} of $q$. Clearly $B_q$ is symmetric in general and alternating in characteristic two. Also if $M$ is a $kG$-module we say $q$ is $G$-invariant if $q(gx)=q(x)$, for all $g\in G$ and $x\in M$. Finally two $G$-invariant quadratic forms $q$ and $q'$ are called \emph{isometric} if there is some $\alpha\in \Aut_{kG}(M)$ so that $q(\alpha(x))=q'(x)$, for all $x\in M$.

Let $n$ be the dimension of $M$ and let $Q\in\Mat_n(k)$. Then $q_Q(x):=x^TQx$, for all $x\in M$, defines a quadratic form on $M$. Its associated bilinear form is $B_q$, which corresponds to the matrix $Q+Q^T$. If $N\in\Mat_n(k)$ is symmetric and hollow, then $q_{Q+N}(x)=x^T(Q+N)x=x^TQx=q_Q(x)$, for all $x\in M$. On the other hand for every quadratic form $q$ on $M$ there is some $Q\in\Mat_n(k)$, such that $q(x)=x^TQx$, for all $x\in M$. Note that $Q$ is not uniquely determined in $\Mat_n(k)$, but within the group of upper-triangular matrices. Now suppose there is another $Q'\in\Mat_n(k)$, such that $q=q_{Q'}$. Then there are symmetric and hollow $N,N'\in\Mat_n(k)$ such that both $Q+N$ and $Q'+N'$ are upper triangular. Since $q_{Q+N}=q=q_{Q'+N'}$ it follows that $Q+N=Q'+N'$, that is, $Q$ and $Q'$ differ by a symmetric and hollow matrix. 

For matrices $Q,Q'\in\Mat_n(k)$ we define $Q\cong Q'$, if $q_Q=q_{Q'}$. Then by the previous paragraph $Q\cong Q'$ if and only if $Q=Q'+N$, for some symmetric and hollow $N\in\Mat_n(k)$. If furthermore $q_Q$ and $q_{Q'}$ are $G$-invariant, then $q_Q$ and $q_{Q'}$ are isometric if and only if $Q\cong \alpha^TQ'\alpha$, for some $\alpha\in \Aut_{kG}(M)$.

For the remainder of this introduction let $k$ have characteristic two. Also let $q$ and $q'$ be quadratic forms on $M$ with corresponding matrices $Q$ and $Q'$ in $\Mat_n(k)$. If $B_q=B_{q'}$, then $Q+Q^T=Q'+Q'^T$, and so $Q+Q'$ is symmetric. In particular $Q'=Q+N+D$, where $N\in\Mat_n(k)$ is symmetric and hollow and $D\in\Mat_n(k)$ is a diagonal matrix. Hence $Q'\cong Q+D$. 

Next let $S\in\Mat_n(k)$ represent a symplectic form on $M$. We define $\widehat{S}\in\Mat_n(k)$ as follows
\begin{equation}\label{S-hat}
\widehat{S}_{i,j}=\begin{cases} S_{i,j}&\text{if $i<j$}\\0&\text{if $i\geq j$}\end{cases}
\end{equation}
As $\widehat{S}+\widehat{S}^T=S$ we see that the quadratic form $q_{\widehat{S}}(x)=x^T\widehat{S} x$, for all $x\in M$, has $S$ as its associated bilinear form. In particular by the previous paragraph it follows that for a quadratic form $q$ on $M$ with corresponding matrix $Q\in\Mat_n(k)$ and associated bilinear form $S$ we have that $Q\cong \widehat{S}+D$, for some diagonal matrix $D\in\Mat_n(k)$. If now in addition $S$ is $G$-invariant the question that remains and which we wish to answer in the case of the Klein-Four Group is, for which $D$, if any, does $\widehat{S}+D$ correspond to a $G$-invariant quadratic form.

Finally let $\EuScript{S}$ denote the isometry class of $S$ and let $Q$ be a $G$-invariant quadratic form with associated bilinear form $S$. Then for any $\alpha\in\Aut_{kG}(M)$ we have that $\alpha^TQ\alpha$ and $Q$ are isometric and $\alpha^TQ\alpha$ has associated bilinear form $\alpha^TS\alpha$ which is in $\EuScript{S}$ but not necessarily equal to $S$. In particular the set of all $G$-invariant quadratic forms with associated bilinear form $S$ is generally not closed under isometry. However every $G$-invariant quadratic form with an associated bilinear form in $\EuScript{S}$ is isometric to a $G$-invariant quadratic form with an associated bilinear form $S$. Hence in order to understand the isometry classes of all $G$-invariant quadratic forms with an associated bilinear form in $\EuScript{S}$ it is sufficient to find a full set of non-isometric representatives for the set of all $G$-invariant quadratic forms with an associated bilinear form $S$.\\

{\bf{Other notation}:} 
Still let $k$ denote a field. A matrix with constant diagonals is called \emph{Toeplitz-matrix} and a matrix with constant anti-diagonals is called \emph{Hankel-matrix}. For $m\times n$-matrices, let $T_i(x)$ be the Toeplitz-matrix where the entries $(s,t)$ with $t-s=i-1$ equal $x\in k$ and all other entries equal zero, and let $H_i(x)$ be the Hankel-matrix where the entries $(s,t)$ with $s+t=i$ equal $x$ and all other entries equal zero. Note that the size of the matrices $T_i(x)$ and $H_i(x)$, though crucial, is not used in the notation but shall always be clear from the context. Furthermore let $\EuScript{T}_{m\times n}(k)$ and $\EuScript{H}_{m\times n}(k)$ denote the set of all $m\times n$-Toeplitz-matrices and all $m\times n$-Hankel-matrices, respectively. Also let $\widetilde{I}_m$ be the $m\times m$-matrix with ones on the anti-diagonal and zeros everywhere else. Then $\widetilde{I}_m\cdot T_i(x)=H_{m+i}(x)$, for all $x\in k$, and $\widetilde{I}_m\cdot \EuScript{T}_{m\times n}(k)=\EuScript{H}_{m\times n}(k)$. 

We call a Toeplitz-matrix that is zero below the main diagonal \emph{upper-triangular} and denote by $\EuScript{UT}_n(k)$ the set of all upper-triangular Toeplitz-matrices in $\EuScript{T}_{n\times n}(k)$. Multiplication in $\EuScript{UT}_n(k)$ is commutative, as for all $x,y\in k$ and integers $i,j\geq 1$ we have $T_i(x)\cdot T_j(y)=T_{i+j-1}(x\cdot y)$, if $i+j-1\leq n$, and $T_i(x)\cdot T_j(y)=0$ otherwise. We call a Hankel-matrix which is zero above the anti-diagonal \emph{lower-triangular} and denote by $\EuScript{LH}_n(k)$ the set of all lower-triangular Hankel-matrices in $\EuScript{N}_{n\times n}(k)$. Note that such matrices are symmetric. For $A\in \EuScript{UT}_n(k)\cup \EuScript{LH}_n(k)$ we have $\rk(A)=0$, if $A=0$ and $\rk(A)=n+1-l$, if $A=\sum_{i=l}^n T_i(A_i)$ or $A=\sum_{i=l}^n H_{n+i}(A_i)$, with $A_l\neq 0$. Next let in addition $B\in \EuScript{UT}_n(k)$ such that $AB\neq 0$. Then we have 
\begin{equation}\label{equation-rank(AB)}
\rk(AB)=\rk(A)+\rk(B)-n.
\end{equation}
We write $E_{s,t}(x)\in\Mat_{m\times n}(k)$ for the single-entry $m\times n$--matrix where all entries are zero except for the $(s,t)$-entry which is $x$. Finally by $\Diag_n(k)$ we denote the set of all diagonal matrices in $\Mat_n(k)$. Also for matrices $A_i\in\Mat_{n_i}(k)$, for $i=1,\ldots,t$, let $\diag(A_1,\ldots,A_t)$ be the matrix in $\Mat_{n_1+\ldots+n_t}(k)$ with $A_1,\ldots,A_t$ on its diagonal. 

\section{Indecomposable symplectic $k(C_2\times C_2)$-modules}\label{Chapter3}
Throughout let $G=C_2\times C_2=\{1,g_1,g_2,g_1g_2\}$ be the Klein Four-Group. Also let $k$ be a perfect field of characteristic two. The Klein-Four Group is of tame representation type and the indecomposable $kG$ modules are described in \cite{Conlon}. They are enumerated using the notation: 
\[\text{$k_G$,\ $kG$,\ $A_n$,\ $B_n$,\ $C_n(f)$ and $C_n(\infty)$, }\]
where $n\geq 1$ is an integer and $f\in k[T]$ is an irreducible polynomial over $k$ in the indeterminate $T$. In \cite{Conlon} the trivial $kG$-module is $A_0=B_0$ and the regular module $kG$ is $D$, but we will instead use the more standard notation $k_G$ and $kG$, respectively.

Recall from the introduction that an indecomposable symplectic $kG$-module is even-dimensional and either indecomposable as a $kG$-module or isomorphic to $N\oplus N^*$, for some indecomposable $kG$-module $N$. As for all integers $n\geq 1$ the modules $A_n$ and $B_n$ are odd-dimensional, $A_n^*\cong B_n$ and $C_n(f)$ and $C_n(\infty)$ are self-dual, we see that the only $kG$-modules that are candidates for indecomposable symplectic modules are:
\[(k_G)^2,\ kG,\ (kG)^2,\ A_n\oplus B_n,\ C_n(f),\ C_n(f)^2,\ C_n(\infty),\ C_n(\infty)^2,\]
where $M^2:=M\oplus M$, for some $kG$-module $M$. During the remainder of this paper we discuss these modules one by one and for each module the results are summarized in the Theorems \ref{Theorem_(k_G)2}, \ref{Theorem_D}, \ref{Theorem_D^2}, \ref{Theorem_AnoplusBn}, \ref{Theorem-C_n(f)}, \ref{Theorem-C_n(f)+C_n(f)*}, \ref{Theorem-C_n(infty)} and \ref{Theorem-C_n(infty)+C_n(infty)*}, respectively.

We also make frequent use of the floor function, that is, $\lfloor x\rfloor$ denotes the largest integer less than or equal to $x\in\R$. Also throughout the paper let $\EuScript{P}$ be a full set of representatives for the distinct cosets of the additive subgroup $\{x^2+x:\ x\in k\}$ in $(k,+)$. In particular $|\EuScript{P}|=2$, if $k$ is finite and $|\EuScript{P}|=1$, if $k$ is algebraically closed. Finally let $\EuScript{IS}(M)$ denote any full set of representatives for the distinct isometry classes of indecomposable symplectic $kG$-modules $(M,S)$ and let $\EuScript{Q}(S)$ denote any full set of representatives for the distinct isometry classes of $G$-invaraint quadratic form with associated bilinear form $S$. 

\subsection{$k_G\oplus k_G$}
As the group action on $k_G\oplus k_G$ is trivial any statement about symplectic indecomposable modules $(k_G\oplus k_G,S)$ is just a statement about $k$-vector spaces. This theory is well understood (e.g. see \cite{Grove}, Theorem 2.10 and Theorem 12.9). We briefly summarize those results and comment on the isometry classes of the quadratic forms.
\begin{Theorem}\label{Theorem_(k_G)2}
We have $\EuScript{IS}(k_G\oplus k_G)=\{S\}$, where $S$ is given on a basis $\{e_1,e_2\}$ for $k_G\oplus k_G$ by $S(e_1,e_2)=1_k$. Next for each $d\in k$ set $q_d(\alpha e_1+\beta e_2)=\alpha^2+\alpha\beta+\beta^2\cdot d$, for all $\alpha,\beta\in k$. Then $\EuScript{Q}(S)=\{q_x:\ x\in \EuScript{P}\}$. 
\end{Theorem}

Proof: Set $V:=k_G\oplus k_G$. By [\cite{Grove}, Theorem 2.10] we know that any two-dimensional symplectic $k$-vector space is a hyperbolic plane. In particular there is one isometry class. Next let $(V,S)$ indecomposable symplectic module and let $\{e_1,e_2\}$ for $V$ such that $S(e_1,e_2)=1_k$. Also let $q$ be a quadratic form on $V$ such that $B_q=S$. Then $q$ is isometric to $q_d$, for some $d\in k$, by [\cite{Grove}, Theorem 12.9], and $q$ is isometric to $q_x$, for a unique $x\in \EuScript{P}$, by [\cite{Grove}, Proposition 13.14]\\
\qed\\\\

\subsection{$kG$ and $kG\oplus kG$}\label{subsection-kG+kG}
Note that \cite{GowWillems2} describes all symplectic and quadratic forms of $kG$, (for any group $G$) but gives no details on isometries. For $x=a\cdot 1_G+b\cdot g_1+c\cdot g_2+d\cdot g_1g_2\in kG$ we write $x=(a,b,c,d)$ and consider the augmentation map $kG\rightarrow k: x\mapsto |x|:=a+b+c+d$, an epimorphism of $k$-algebras. Note that $x^2=|x|^2\cdot 1_G$. Hence $x$ is a unit in $kG$ if $|x|\neq 0_k$, and $x$ is nilpotent otherwise. Also recall that $\End_{kG}(kG)\cong kG$, via $\varphi\leftrightarrow \varphi(1_G)$. Next let $B$ be the $G$-invariant symmetric form on $kG$ where $B(g,h)=\delta_{g,h}$, for all $g,h\in G$, is extended $k$-linearly on $kG$. Now for $\varphi\in kG$ set $B_{\varphi}(x,y):=B(\varphi x,y)$, for all $x,y\in kG$. Then $\{B_{\varphi}:\varphi\in kG\}$ gives the space of all $G$-invariant bilinear forms on $kG$. Note that all $B_{\varphi}$ are symmetric. Furthermore $B_{\varphi}$ is alternating if and only if $\varphi=(0,b,c,d)\in kG$, and $B_{\varphi}$ is non-degenerate if and only if $\varphi$ is a unit in $kG$. Overall it follows that the space of $G$-invariant symplectic forms on $kG$ is given by $\{B_{\varphi}: \varphi=(0,b,c,d)\in kG, |\varphi|\neq 0_k\}$.

Next let $B_{\varphi}$ be a $G$-invariant symplectic forms on $kG$, where $\varphi=(0,\beta,\gamma,\delta)$ is a unit in $kG$. Then $\{q_{r,\varphi}: r\in k\}$ gives all $G$-invariant quadratic form with associated bilinear form $B_{\varphi}$, where 
\begin{equation}\label{equation-qrvarphi}
q_{r,\varphi}(a,b,c,d):=r\cdot (a+b+c+d)^2+(ab+cd)\beta+(ac+bd)\gamma+(ad+bc)\delta,
\end{equation}
for all $(a,b,c,d)\in kG$ and $r\in k$.

\begin{Theorem}\label{Theorem_D}
We have $\EuScript{IS}(kG)=\{B_{\varphi}: \varphi=(0,b,c,1+b+c)\in kG\}$ and $\EuScript{Q}(B_{\varphi})=\{q_{0_k,\varphi}\}$, for any unit $\varphi=(0,b,c,d)$ in $kG$. In particular if $k$ is finite, then $|\EuScript{IS}(kG)|=|k|^2$ and there is a single isometry class of $G$-invariant quadratic forms with associated bilinear form $B_{\varphi}$.
\end{Theorem}

Proof: Let $\varphi=(0,b,c,d)$ and $\varphi'=(0,b',c',d')$ be units in $kG$. Then $B_{\varphi}$ and $B_{\varphi'}$ are isometric if and only if there exists $\psi\in \Aut_{kG}(kG)$ so that $B_{\varphi}(\psi(x),\psi(y))=B_{\varphi'}(x,y)$, for all $x,y\in kG$. Note that $\psi(1_G)$ is a unit, $[\psi(1_G)]^2=|\psi(1_G)|^2\cdot 1_G$ and $B_{\varphi}(\psi(x),\psi(y))=B_{\varphi}(\psi(1_G)x,\psi(1_G)y)=B_{\varphi [\psi(1_G)]^2}(x,y)$. Thus $B_{\varphi}$ and $B_{\varphi'}$ are isometric if and only if $\varphi$ and $\varphi'$ are non-zero multiples of each other. Thus for every indecomposable symplectic module $(kG,S)$ there is a unique $\varphi=(0,b,c,1+b+c)\in kG$ such that $S$ is isometric to $B_{\varphi}$.

Next let $\varphi=(0,b,c,d)$ be a unit in $kG$ and let $r\in k$ be given. Then $|\varphi|\neq 0$. Also let $\zeta=(1_G,1_G,1_G,1_G)\in kG$. Then there is $\psi\in\Aut_{kG}(kG)$ so that $\psi(1_G)=1_G+r|\varphi|^{-1}\cdot \zeta\in kG$. Now for all $\alpha\in kG$ we have $\psi(\alpha)=\alpha+t\zeta$, where $t=r|\varphi|^{-1}\cdot |\alpha|$. One checks that $B_{\varphi}(\alpha,\zeta)=|\varphi|\cdot |\alpha|$. Thus
\begin{align*}
q_{0_k,\varphi}(\psi(\alpha))&=q_{0_k,\varphi}(\alpha+t\zeta)=q_{0_k,\varphi}(\alpha)+t^2\cdot \underbrace{q_{0_k,\varphi}(\zeta)}_{=0_k}+ \underbrace{t\cdot B_{\varphi}(\alpha,\zeta)}_{=r\cdot |\alpha|^2}=q_{r,\varphi}(\alpha).
\end{align*}
\qed\\\\

Let us now turn to $kG\oplus kG$. The augmentation map is a homomorphism with kernel $J(kG)$. Hence $kG/J(kG)\cong k$ and so $\End_{kG}(kG)/J(\End_{kG}(kG))\cong k$. Next observe that with respect to the basis $\{1_G,g_1,g_2,g_1g_2\}$ the element $x=(a,b,c,d)\in kG$ corresponds to the matrix
\begin{equation}\label{equation-xAsMatrix}
\hat{x}=a\cdot I_4+b\cdot \begin{pmatrix} \widetilde{I_2}&0\\0&\widetilde{I_2}\end{pmatrix}+c\cdot \begin{pmatrix} 0&I_2\\I_2&0\end{pmatrix}+d\cdot \begin{pmatrix} 0&\widetilde{I_2}\\\widetilde{I_2}&0\end{pmatrix}\in \Mat_4(k).
\end{equation}
In this sense $kG$ is injectively embedded in $\Mat_4(k)$ and henceforth we identify $x\in kG$ with its corresponding matrix $\hat{x}\in\Mat_4(k)$. Also note that $\hat{x}$ is self-transpose and its top row equals the vector $(a,b,c,d)$. 

For $\alpha,\beta,\gamma\in kG\subseteq \Mat_4(k)$ we define $S(\alpha,\beta,\gamma):=\begin{pmatrix} \alpha&\beta\\\beta&\gamma\end{pmatrix}\in\Mat_8(k)$. Furthermore set $\EuScript{A}:=\{g_1+\lambda g_2+ (1+\lambda) g_1g_2:\ \lambda\in k\}$ and 
\[\EuScript{S}:=\{S(0,1_G,0)\}\cup \{S(\alpha,1_G,0): \alpha\in\EuScript{A}\}\cup\{S(g_2+g_1g_2,1_G, \mu(g_1+g_1g_2)): \mu\in k\}.\]
In the following we represent all forms on $kG\oplus kG$ with respect to the basis $\EuScript{B}:=\{(1_G,0), (g_1,0), (g_2,0), (g_1g_2,0),(0,1_G), (0,g_1), (0,g_2),(0,g_1g_2)\}$.

\begin{Theorem}\label{Theorem_D^2}
(a) {\bf{symplectic forms}}: We have $\EuScript{IS}(kG\oplus kG)=\EuScript{S}$. In particular if $k$ is finite, then $|\EuScript{IS}(kG\oplus kG)|=2|k|+1$.\\
\hfill\\
(b) {\bf{quadratic forms}}: Next let $(kG\oplus kG,S)$ be an indecomposable symplectic $kG$-module. Then the $G$-invariant quadratic forms with associated bilinear form $S$ are precisely given by $Q_{x,y}:=\widehat{S}+\diag(x\cdot I_4,y\cdot I_4)\in\Mat_8(k)$, for all $x,y\in k$. Finally $\EuScript{Q}(S)=\{Q_{0,0}\}$.
\end{Theorem}

Proof: (a) {\bf{1. Notation}}: Recall the bijection between $kG$ and the set of $G$-invariant bilinear forms on $kG$, via $\varphi\leftrightarrow B_{\varphi}$, as explained in the first paragraph of this section. Next by Remark \ref{Remark-symplectic_forms_on_M^2} the indecomposable symplectic $kG$-modules $(kG\oplus kG,S)$ are precisely given by $S=S(\alpha,\beta,\gamma)$, where $\alpha,\beta,\gamma\in kG$ so that $\alpha$ and $\gamma$ correspond to alternating degenerate $G$-invariant forms $B_{\alpha}$ and $B_{\gamma}$ on $kG$ and $\beta$ corresponds to a $G$-invariant form $B_{\beta}$ on $kG$. Thus $\alpha=(0,\alpha_2,\alpha_3,\alpha_2+\alpha_3)$, $\gamma=(0,\gamma_2,\gamma_3,\gamma_2+\gamma_3)$ and $\beta=(\beta_1,\beta_2,\beta_3,\beta_4)$. Furthermore note that since $S$ is non-degenerate it follows that $\alpha\gamma+\beta^2$ is a unit in $kG$. As $\alpha$ and $\gamma$ are non-units this forces that $\beta$ is a unit. In the following let $S=S(\alpha,\beta,\gamma)$ be such a form.\\
\hfill\\
\noindent {\bf{2. A first reduction}}: Next let $M\in\Aut_{kG}(kG\oplus kG)$. Then $M=\begin{pmatrix} a&b\\c&d\end{pmatrix}\in\Mat_8(k)$, where $a,b,c,d\in\End_{kG}(kG)\cong kG$ and $ad+bc$ is a unit. Now
\begin{align*}
M^TSM&=\begin{pmatrix} a&c\\b&d\end{pmatrix}\cdot\begin{pmatrix} \alpha&\beta\\\beta&\gamma\end{pmatrix}\cdot\begin{pmatrix} a&b\\c&d\end{pmatrix}\\&=\begin{pmatrix} a^2\alpha+c^2\gamma&ab\alpha+cd\gamma+(ad+bc)\beta\\ ab\alpha+cd\gamma+(ad+bc)\beta&b^2\alpha+d^2\gamma\end{pmatrix}
\end{align*}
First note that multiplication by $x^2$, for some $x\in kG$ is equivalent to multiplication by $|x|^2\in k$. Furthermore since $k$ is perfect, for every $x'\in k$ there is $x_1\in k$ such that $x_1^2=x'$ and thus $(x_1\cdot 1_G)^2=x'\cdot 1_G$. Also given a symplectic $G$-invariant form $S(\alpha,\beta,\gamma)$ let $a=1_G$, $b=c=0$ and $d=\beta^{-1}$. Then $M$ is invertible and thus $S(\alpha,\beta,\gamma)$ is isometric to $S(\alpha,1_G,\mu\gamma)$, for some $\mu\in k^*$. (Here $\mu$ satisfies $d^2=\mu\cdot 1_G$).\\
\hfill\\
\noindent {\bf{3. A relation between isometric forms}}: Clearly the isometry class of $S(0,1_G,0)$ is the set $\{S(0,\beta',0):\ \text{$\beta'\in kG$ a unit}\}$. In particular $S(0,1_G,0)$ is not isometric to any other element in $\EuScript{S}$. From now on we assume that $\alpha$ and $\gamma$ are not both zero. Also $S(\alpha,\beta,\gamma)$ is isometric to $S(\alpha',\beta',\gamma')$, for some unit $\beta'\in kG$, if and only if $\alpha'=a^2\alpha+c^2\gamma$ and $\gamma'=b^2\alpha+d^2\gamma$. Let $\alpha'=(0,\alpha_2',\alpha_3',\alpha_2'+\alpha_3')$ and $\gamma'=(0,\gamma_2',\gamma_3',\gamma_2'+\gamma_3')$. Then $S(\alpha,\beta,\gamma)$ is isometric to $S(\alpha',\beta',\gamma')$, for some unit $\beta'\in kG$, if and only if 
\begin{equation}\label{Equation-Isometry}
\begin{pmatrix} \alpha_2&\gamma_2\\\alpha_3&\gamma_3 \end{pmatrix}\cdot \begin{pmatrix} a'&b'\\c'&d'\end{pmatrix}=\begin{pmatrix} \alpha_2'&\gamma_2'\\\alpha_3'&\gamma_3' \end{pmatrix},
\end{equation}
for some $a',b',c',d'\in k$ with $a'd'+b'c'\neq 0_k$. In the following we distinguish between $\alpha_2\gamma_3+\alpha_3\gamma_2$ being different from and equal to zero.\\
\hfill\\
\noindent {\bf{4. The case $\alpha_2\gamma_3+\alpha_3\gamma_2\neq 0_k$}}: Then for $\alpha_2'=0_k$, $\alpha_3'=1_k$, $\gamma_2'=1_k$ and $\gamma_3'=0_k$ there are $a',b',c',d'\in k$ with $a'd'+b'c'\neq 0$ such that (\ref{Equation-Isometry}) holds. In particular $S(\alpha,\beta,\gamma)$ is isometric to $S(g_2+g_1g_2,\beta',g_1+g_1g_2)$, for some unit $\beta'\in kG$. Hence by our second step $S(\alpha,\beta,\gamma)$ is isometric to $S(g_2+g_1g_2,1_G,\mu(g_1+g_1g_2))$, for some $\mu\in k^*$.

For $\mu\in k^*$ let $S_{\mu}=S(g_2+g_1g_2,1_g,\mu(g_1+g_1g_2))$. Then $S_{\mu}$ cannot be isometric to some $S(\alpha',\beta',0)$ as otherwise taking the determinant in (\ref{Equation-Isometry}) leads to a contradiction. Now assume that $S_{\mu}$ and $S_{\mu'}$ are isometric for $\mu,\mu'\in k^*$. Then 
\begin{align*}
g_2+g_1g_2&=a^2\cdot (g_2+g_1g_2)+c^2\cdot (g_1+g_1g_2)\\
\mu'(g_1+g_1g_2)&=b^2\cdot (g_2+g_1g_2)+d^2\cdot \mu(g_1+g_1g_2)\\
1_G&=ab\cdot (g_2+g_1g_2)+cd\cdot (g_1+g_1g_2)+(ad+bc)\cdot 1_G.
\end{align*}
The first two equations imply $a^2=1_G$, $c^2=0$ and $b^2=0$, $d^2\mu=\mu'$, respectively. Squaring the third equation gives $1_G=a^2d^2\cdot 1_G=d^2$. Thus $\mu=\mu'$ ensues. Overall we can conclude that $S_{\mu}$ is not isometric to any other element in $\EuScript{S}$.\\
\hfill\\
\noindent {\bf{5. The case $\alpha_2\gamma_3+\alpha_3\gamma_2=0_k$}}: Observe that $S(\alpha,\beta,\gamma)$ is isometric to $S(\gamma,\beta,\alpha)$ with $a=d=0$ and $b=c=1_G$. Hence we can assume that $\alpha\neq 0$. Now there is some $x\in k$ so that $x\alpha\in \EuScript{A}\cup \{g_2+g_1g_2\}$ and there is some $y\in k^*$ so that $\gamma=y\alpha$. Next let $a',b'\in k$ such that $a'^2=x$ and $b'^2=y$. Then with $a=a'\cdot 1_G$, $b=b'\cdot 1_G$, $c=0_k\cdot 1_G$ and $d=1_k\cdot 1_G$ we get that $S(\alpha,\beta,\gamma)$ is isometric to $S(a'^2\alpha,\beta',0)$, for some unit $\beta'\in kG$. Again by step 2 we get that $S(\alpha,\beta,\gamma)$ is isometric to $S(\alpha,1_G,0)$, for some $\alpha\in\EuScript{A}\cup \{g_2+g_1g_2\}$. Finally note that if $S(\alpha,1_G,0)$ and $S(\alpha',1_G,0)$ are isometric for $\alpha,\alpha'\in \EuScript{A}\cup \{g_2+g_1g_2\}$, then there is some $a\in kG$ such that $a^2\alpha=\alpha'$, which forces $\alpha=\alpha'$. \\

(b) {\bf{1. $G$-invariant quadratic forms}}: Let $q$ be a $G$-invariant quadratic form with associated bilinear form $S$. Then $q(g,0)=q(1_G,0)$ and $q(0,g)=q(0,1_G)$, for all $g\in G$. Thus $q$ corresponds to the matrix $Q_{x,y}:=\widehat{S}+\diag(x\cdot I_4,y\cdot I_4)\in\Mat_8(k)$, where $x=q(1_G,0)$ and $y=q(0,1_G)$. On the other hand suppose that $q$ corresponds to $Q_{x,y}$, for some $x,y\in k$. Then $q$ is a quadratic form with associated bilinear form $S$ and $q=q_{x,y}$, where
\[q_{x,y}(a,b)=q_{x,\alpha}(a)+q_{y,\gamma}(b)+B_{\beta}(a,b),\]
for all $(a,b)\in kG\oplus kG$, where $q_{x,\alpha}$ and $q_{y,\gamma}$ are the $G$-invariant quadratic forms given by (\ref{equation-qrvarphi}). In particular notice that $q_{x,y}$ is $G$-invariant. Overall $\{Q_{x,y}: x,y\in k\}$ gives a complete set of all $G$-invariant quadratic forms with associated bilinear form $S$.\\
\hfill\\
\noindent {\bf{2. isometry class}}: Let $\psi_{x,y}(a, b)=(a+x\cdot \zeta\cdot b,b+y\cdot \zeta\cdot a)$, where $\zeta=(1_G,1_G,1_G,1_G)\in kG$. Then $\psi_{x,y}=\begin{pmatrix} 1_G& x\cdot \zeta\\ y\cdot\zeta&1_G\end{pmatrix}\in \Aut_{kG}(kG\oplus kG)$. Now 
\begin{align*}
q_{0,0}(\psi_{x,y}&(a, b))=q_{0,\alpha}(a+x\cdot \zeta\cdot b)+q_{0,\gamma}(b+y\cdot \zeta\cdot a)+B_{\beta}(a+x\cdot \zeta\cdot b,b+y\cdot \zeta\cdot a)\\
&=q_{0,\alpha}(a+x\cdot \zeta\cdot b)+q_{0,\gamma}(b+y\cdot \zeta\cdot a)+B_{\beta}(a,b)\\&\ +y\cdot B_{\beta}(a,\zeta\cdot a)+x\cdot B_{\beta}(\zeta\cdot b,b)+(xy)\cdot B_{\beta}(\zeta\cdot b,\zeta\cdot a)
\end{align*}
Now one checks that $q_{0,\alpha}(a+x\cdot \zeta\cdot b)=q_{0,\alpha}(a)$, $q_{0,\gamma}(b+y\cdot \zeta\cdot a)=q_{0,\gamma}(b)$, $B_{\beta}(a,\zeta\cdot a)=|a|^2\cdot |\beta|\cdot 1_G$, $B_{\beta}(\zeta\cdot b,b)=|b|^2\cdot |\beta|\cdot 1_G$ and $B_{\beta}(\zeta\cdot b,\zeta\cdot a)=|a|\cdot |b|\cdot |\beta|\cdot B(\zeta,\zeta)=0$. Hence
\begin{align*}
q_{0,0}(\psi_{x,y}&(a, b))=q_{0,\alpha}(a)+x\cdot |a|^2\cdot |\beta|\cdot 1_G+q_{0,\gamma}(b)+y\cdot |b|^2\cdot |\beta|\cdot 1_G+B_{\beta}(a,b)\\
&=q_{x\cdot |\beta|,\alpha}(a)+q_{y\cdot |\beta|,\gamma}(b)+B_{\beta}(a,b)=q_{x\cdot |\beta|,y\cdot |\beta|}(a, b)
\end{align*}
As $|\beta|\neq 0$ it follows that $q_{0,0}$ is isometric to $q_{x,y}$, for all $x,y\in k$.\\
\qed\\\\

\subsection{$A_n\oplus B_n$}
Recall Toeplitz-matrices, Hankel-matrices and single-entry matrices as described near the end of chapter \ref{Chapter2}. In the following we describe the first two families of indecomposable $kG$-modules as given by \cite{Conlon}. Let $n\geq 1$ be an integer and let $A_n=B_n=k^{2n+1}$ as $k$-vector spaces with the basis $\{e_1,\ldots,e_{2n+1}\}$ of unit coordinate column vectors. The action of $G$ on $A_n$ is given by 
\[g_1=\begin{pmatrix} I_n& T_2(1)\\ 0&I_{n+1}\end{pmatrix}\ \text{and}\ g_2=\begin{pmatrix} I_n& T_1(1)\\ 0&I_{n+1}\end{pmatrix}.\]
Here $T_1(1),T_2(1)\in\EuScript{T}_{n\times (n+1)}(k)$. The $G$-action on $B_n$ is given by 
\[g_1=\begin{pmatrix} I_{n+1}& T_1(1)\\ 0&I_n\end{pmatrix}\ \text{and}\ g_2=\begin{pmatrix} I_{n+1}& T_0(1)\\ 0&I_n\end{pmatrix}.\]
Here $T_0(1),T_1(1)\in\EuScript{T}_{(n+1)\times n}(k)$. Clearly $A_n$ and $B_n$ are odd-dimensional and thus have no symplectic forms. 
\begin{Lemma}
For all integers $n\geq 1$ we have $A_n^*\cong B_n$.
\end{Lemma}

Proof: Let $\alpha:G\rightarrow \GL_{2n+1}(k)$ and $\beta:G\rightarrow \GL_{2n+1}(k)$ denote the indecomposable representations $A_n$ and $B_n$ of $G$, respectively. Recall that the dual representation $\alpha^*$ is given by $\alpha^*(g)=\alpha(g^{-1})^T$, for all $g\in G$. Hence with respect to the dual basis $\{e_1^*,\ldots,e_{2n+1}^*\}$ the $G$-action on $A_n^*$ is given by
\[\alpha^*(g_1)=\begin{pmatrix} I_n& 0\\ T_0(1)&I_{n+1}\end{pmatrix}\ \text{and}\ \alpha^*(g_2)=\begin{pmatrix} I_n& 0\\ T_1(1)&I_{n+1}\end{pmatrix}.\] 
As $\widetilde{I}_{2n+1} \cdot \alpha^*(g_1)\cdot \widetilde{I}_{2n+1} =\beta(g_1)$ and $\widetilde{I}_{2n+1} \cdot \alpha^*(g_2)\cdot \widetilde{I}_{2n+1} =\beta(g_2)$, we get $A_n^*\cong B_n$.\\
\qed\\\\

\newcommand\mydots{\hbox to 1em{.\hss.\hss.}}
In the following let $V_n:=A_n\oplus B_n$, for $n\geq 0$. For reasons of presentation we work with respect to the basis 
\[\EuScript{B}:=\{(e_1\oplus 0,\mydots,e_n\oplus 0,0\oplus e_1,\mydots,0\oplus e_{n+1}, e_{n+1}\oplus 0,\mydots,e_{2n+1}\oplus 0,0\oplus e_{n+2},\mydots,0\oplus e_{2n+1}\}.\]
Then the $G$-action on $V_n$ is given by
\[g_1=\begin{pmatrix} I_{2n+1}&T_2(1)\\0&I_{2n+1}\end{pmatrix}\text{ and }g_2=\begin{pmatrix} I_{2n+1}&I_{2n+1}+E_{n+1,n+1}(1)\\0&I_{2n+1}\end{pmatrix}.\]

For any $A\in\Mat_{n\times m}(k)$ we define $A^{\tau}:=\widetilde{I}_m A^T \widetilde{I}_n$. This variation of a transpose turns the first (second etc.) column of $A$ into the last (second last etc.) row of $A^{\tau}$ but in reversed order. In the following we use $\star$ to denote an arbitrary matrix over $k$ of the obvious size.
\begin{Lemma}
\[\End_{kG}(V_n)=\left\{\begin{pmatrix}\alpha I_n&\Sigma&\star&\star\\&\beta I_{n+1}&\star&\star\\&&\alpha I_{n+1}&\Sigma^{\tau}\\&&&\beta I_n\end{pmatrix}:\ \begin{matrix}\alpha,\beta\in k,\\ \Sigma\in \EuScript{T}_{n\times n+1}(k)\end{matrix}\right\}.\]
\end{Lemma}

Proof: Let $M=\begin{pmatrix} A&B\\C&D\end{pmatrix}\in \Mat_{4n+2}(k)$, where $A,B,C,D\in\Mat_{2n+1}(k)$. Then $M\in\End_{kG}(V_n)$ if and only if $Mg_1=g_1M$ and $Mg_2=g_2M$, that is, if and only if 
\begin{equation*}
\begin{split}
(a)\quad &A\cdot T_2(1)=T_2(1)\cdot D\\
(b)\quad &C\cdot T_2(1)=0=T_2(1)\cdot C\\
(c)\quad &A+A\cdot E_{n+1,n+1}(1)=D+E_{n+1,n+1}(1)\cdot D\\
(d)\quad &C+C\cdot E_{n+1,n+1}(1)=0=C+E_{n+1,n+1}(1)\cdot C
\end{split}
\end{equation*}
First note that (b) and (d) hold if and only if $C=0$. Next (a) holds precisely if for all $k\in\{1,\ldots,2n\}$ and $l\in\{2,\ldots,2n+1\}$ we have 
\begin{equation}\label{Equation-A_k,l-1=D_k+1,l}
\text{$A_{k,l-1}=D_{k+1,l}$ and $A_{2n+1,k}=0=D_{l,1}$.}
\end{equation}
Finally (c) holds if and only if for all $k,l\in\{1,\ldots,n,n+2,\ldots,2n+1\}$ we have 
\begin{equation}\label{Equation-A_k,l=D_k,l}
\text{$A_{k,l}=D_{k,l}$ and $A_{n+1,l}=0=D_{k,n+1}$.}
\end{equation}
Hence the matrices described in the statement indeed commute with $g_1$ and $g_2$ and thus belong to $\End_{kG}(V_n)$.

Next suppose that $C=0$ and that $M$ satisfies (\ref{Equation-A_k,l-1=D_k+1,l}) and (\ref{Equation-A_k,l=D_k,l}). By (\ref{Equation-A_k,l=D_k,l}) we have $A=\begin{pmatrix}X_1&a_1&X_2\\0_{1\times n}&a_2&0_{1\times n}\\X_3&a_3&X_4\end{pmatrix}$ and $D=\begin{pmatrix}X_1&0_{n\times 1}&X_2\\d_1&d_2&d_3\\X_3&0_{n\times 1}&X_4\end{pmatrix}$, where $X_i\in \Mat_n(k)$, for $i=1,2,3,4$, and $a_1,a_3\in\Mat_{n\times 1}(k)$, $d_1,d_3\in\Mat_{1\times n}(k)$ and $a_2,d_2\in k$. Let $k,l\in\{1,\ldots,2n+1\}$. Then $A_{k,l}=D_{k+1,l+1}$, unless $k=2n+1$ or $l=2n+1$, by (\ref{Equation-A_k,l-1=D_k+1,l}). Also $A_{k,l}=D_{k,l}$, unless $k=n+1$ or $l=n+1$, by (\ref{Equation-A_k,l=D_k,l}). In particular $A_{k,l}=A_{k+1,l+1}$, unless $k,l\in\{n,2n+1\}$ and $D_{k,l}=D_{k+1,l+1}$, unless $k,l\in\{n+1,2n+1\}$. Therefore the following matrices are Toeplitz-matrices
\[\begin{pmatrix}X_1&0_{n\times 1}\\d_1&d_2\end{pmatrix},\ \begin{pmatrix}a_1&X_2\\\star&d_3\end{pmatrix},\ \begin{pmatrix}0_{1\times n}\\X_3\end{pmatrix},\ \begin{pmatrix}a_2&0_{1\times n}\\a_3&X_4\end{pmatrix}.\]
Above we have seen that the first column of $D$ is zero except possibly for the first entry, which we call $\alpha$, and the last row of $A$ is zero except possibly for the last entry, which we call $\beta$. Thus it follows that $X_1=\alpha I_n$, $d_1=0_{1\times n}$, $d_2=\alpha$, $X_3=0_{n\times n}$, $a_2=\beta$, $a_3=0_{n\times 1}$, $X_4=\beta I_n$ and $d_3=a_2^{\tau}$. Finally set $\Sigma:=\begin{pmatrix}a_1&X_2\end{pmatrix}\in\Mat_{n\times (n+1)}(k)$. Then $\Sigma\in \EuScript{T}_{n\times n+1}(k)$ and indeed $\begin{pmatrix}X_2\\d_3\end{pmatrix}=\Sigma^{\tau}$. This completes the proof.\\
\qed\\\\

\begin{Lemma}\label{Lemma-Ginvariant_forms_An+Bn}
The $G$-invariant forms on $V_n$ are precisely the matrices
\[F=\begin{pmatrix}&\begin{matrix}& \beta\widetilde{I}_n\\\alpha\widetilde{I}_{n+1}&\Omega\end{matrix}\\\begin{matrix}&\beta\widetilde{I}_{n+1}\\ \alpha\widetilde{I}_n&\Omega^T\end{matrix}&D\end{pmatrix},\]
where $\alpha,\beta\in k$, $\Omega\in \EuScript{H}_{n+1\times n}(k)$ and $D\in\Mat_{2n+1}(k)$. Furthermore $F$ corresponds to a symplectic form if and only if $\alpha=\beta\neq 0$ and $D$ is hollow and symmetric.
\end{Lemma}

Proof: One checks easily that $g_1^T\cdot \widetilde{I}_{4n+2}\cdot g_1=\widetilde{I}_{4n+2}$ and $g_2^T\cdot \widetilde{I}_{4n+2}\cdot g_2=\widetilde{I}_{4n+2}$. Hence $\widetilde{I}_{4n+2}$ is a $G$-invariant form and thus $\{\widetilde{I}_{4n+2}\cdot A:\ A\in\End_{kG}(V_n)\}$ gives the set of all $G$-invariant forms. Also note that $\Omega:=\widetilde{I}_{n+1} \Sigma^{\tau}=\Sigma^T\widetilde{I}_n\in \EuScript{H}_{n+1\times n}(k)$, for $\Sigma\in\EuScript{T}_{n\times n+1}(k)$, and $(\widetilde{I}_n \Sigma)^T=\widetilde{I}_{n+1} \Sigma^{\tau}$. Hence the statement follows.\\
\qed\\\\

We define two subsets of $\EuScript{H}_{(n+1)\times n}(k)$. First let $\EuScript{F}$ denote the set of all $(n+1)\times n$--Hankel-matrices where the entry of the top-most non-zero anti-diagonal is $1_k$. If, beginning at one, we count the anti-diagonals of an element in $\EuScript{H}_{(n+1)\times n}(k)$ starting in the top left corner and proceeding to the bottom right corner, then we count $2n$ anti-diagonals. Now by $\EuScript{F}'$ we denote those elements in $\EuScript{F}$ where the entry for each even-numbered anti-diagonal equals the entry of the ensuing odd-numbered anti-diagonal, that is, $\EuScript{F}':=\{X\in\EuScript{F}: X_{i,i}=X_{i,i-1}, \text{for all $i=2,\ldots,2n$}\}$. Note that $0\in\EuScript{F}'\subseteq \EuScript{F}$. Also observe that in $\EuScript{F}'$ the top most non-zero anti-diagonal, if it exists, is either the first or an even numbered one.

For $\Omega\in \EuScript{H}_{n+1\times n}(k)$ let $S(\Omega)$ correspond to the symplectic form as described in Lemma \ref{Lemma-Ginvariant_forms_An+Bn} with $\alpha=\beta=1_k$ and $D=0$. 
We can now state the main theorem of this section.
\begin{Theorem}\label{Theorem_AnoplusBn} 
(a) {\bf{symplectic forms:}} We have $\EuScript{IS}(A_n\oplus B_n)=\{S(\Omega): \Omega\in\EuScript{F}\}$. In particular if $k$ is finite then $|\EuScript{IS}(A_n\oplus B_n)|=|k|^{2n-1}+\ldots+|k|+2$.\\
\hfill\\
(b) {\bf{quadratic forms:}} For $\Omega\in\EuScript{F}$ there is a $G$-invariant quadratic form with associated bilinear form $S(\Omega)$ if and only if $\Omega\in\EuScript{F}'$. In this case all such quadratic forms are given by the set $\{Q_{\Omega}(D'): D'\in\Diag_{2n+1}(k)\}$, where $Q_{\Omega}(D'):=\widehat{S(\Omega)}+\diag(D,D')$, and $D\in\Mat_{2n+1}(k)$ is the diagonal matrix such that
\begin{equation}\label{EquationForD1}
D_{i,i}=\begin{cases} 0,&\text{if $i=1,\ldots,n$}\\ \Omega_{i-n,i-n},&\text{if $i=n+1,\ldots,2n$}\\\Omega_{n+1,n},&\text{if $i=2n+1$.}\end{cases}
\end{equation}
If $\Omega=0$, then $\EuScript{Q}(S(\Omega))=\{Q_{\Omega}(0)\}$. Otherwise $\EuScript{Q}(S(\Omega))=\{Q_{\Omega}(E_{t,t}((\eta_{t,t})^{-1} x)): x\in\EuScript{P}\}$, where $\eta:=\widetilde{I}_{2n+1}D\widetilde{I}_{2n+1}$ and $t\in\{1,\ldots,n+1\}$ is minimal such that $\eta_{t,t}\neq 0$.
\end{Theorem}

Proof: (a) {\bf{1. Notation}}: By Lemma \ref{Lemma-Ginvariant_forms_An+Bn} the indecomposable symplectic $kG$-modules $(V_n,S)$ are precisely given by $S(A,D):=\begin{pmatrix} 0&A\\A^T&D \end{pmatrix}$, where $D\in\Mat_{2n+1}(k)$ is symmetric and hollow and $A=A(x,\Omega'):=\begin{pmatrix} 0&x\widetilde{I}_n\\x\widetilde{I}_{n+1}&\Omega'\end{pmatrix}\in\Mat_{2n+1}(k)$, with $x\in k^*$ and $\Omega'\in \EuScript{H}_{n+1\times n}$. In the following let $S=S(A(x,\Omega),D)$ be such a form. Also let $M=\begin{pmatrix} M_1&M_2\\0&M_4 \end{pmatrix}\in\Aut_{kG}(V_n)$, where $M_1=\begin{pmatrix} \alpha I_n&\Sigma\\0&\beta I_{n+1}\end{pmatrix}\in\Mat_{2n+1}(k)$ and $M_4=\begin{pmatrix} \alpha I_{n+1}&\Sigma^{\tau}\\0&\beta I_n\end{pmatrix}\in\Mat_{2n+1}(k)$, for $\alpha,\beta\in k^*$ and $\Sigma\in\EuScript{T}_{n\times n+1}$, and $M_2\in \Mat_{2n+1}(k)$. Then
\begin{equation}\label{Equation-M^TSM}
\begin{split}M^TSM&=\begin{pmatrix} M_1^T&0\\M_2^T&M_4^T \end{pmatrix}\cdot \begin{pmatrix} 0&A\\A^T&D \end{pmatrix}\cdot \begin{pmatrix} M_1&M_2\\0&M_4 \end{pmatrix}\\&=\begin{pmatrix} 0&M_1^TAM_4\\M_4^TA^TM_1&M_2^TAM_4+M_4^TA^TM_2+M_4^TDM_4 \end{pmatrix}
\end{split}
\end{equation}
\hfill\\
\noindent{\bf{2. We show that $S$ is isometric to $S(\Omega)$, for a unique $\Omega\in\EuScript{F}$:}} As $D$ is symmetric and hollow there is some $Z\in\Mat_{2n+1}(k)$ such that $D=Z+Z^T$. Hence if we set $M_1=M_4=I_{2n+1}$ and $M_2=A^{-T}Z$, then it follows from (\ref{Equation-M^TSM}) that $S$ is isometric to $S(A,0)$. Thus from now on we may assume that $D=0$. Also we deduce from (\ref{Equation-M^TSM}) that $S(A,0)$ is isometric to some $S(B,0)$, if and only if $B=M_1^TAM_4$, for $M_1$ and $M_4$ as described in step one. We have
\[M_1^TAM_4=\begin{pmatrix} 0&\alpha\beta x\widetilde{I}_n\\ \alpha\beta x\widetilde{I}_{n+1}&\beta x (\Sigma^T\widetilde{I}_n+\widetilde{I}_{n+1}\Sigma^{\tau})+\beta^2\Omega'\end{pmatrix}=\begin{pmatrix} 0&\alpha\beta x\widetilde{I}_n\\ \alpha\beta x\widetilde{I}_{n+1}&\beta^2\Omega'\end{pmatrix},\]
as $\widetilde{I}_{n+1}\Sigma^{\tau}=\Sigma^T\widetilde{I}_n$. Note that there is some $\beta\in k^*$, such that $\beta^2\Omega'\in \EuScript{F}$. Set $\alpha:=(\beta x)^{-1}\in k^*$. Then $M_1^TAM_4=A(1_k,\Omega)$,
for some $\Omega\in \EuScript{F}$. Furthermore note that $\Omega\in \EuScript{F}$ is unique as $\beta^2\Omega'\in \EuScript{F}$ and $(\beta')^2\Omega'\in \EuScript{F}$, for different $\beta,\beta'\in k^*$, implies that $\Omega'=0$, and thus $\beta^2\Omega'=(\beta')^2\Omega'$.\\
\hfill\\
\noindent (b) {\bf{1. Notation}}: Let $S=S(\Omega)$. Then $\widehat{S}=\begin{pmatrix} 0&A\\ 0&0\end{pmatrix}$, where $A=A(1_k, \Omega)\in \Mat_{2n+1}(k)$. Then the quadratic forms with associated bilinear form $S$ are precisely the matrices $\begin{pmatrix} D_1&A\\ 0&D_2\end{pmatrix}$, where $D_1,D_2\in\Diag_{2n+1}(k)$. Let $Q$ be such a matrix.\\
\hfill\\
\noindent {\bf{2. $G$-invariant quadratic forms}}: Set $P_1:=T_2(1)\in\Mat_{2n+1}(k)$ and $P_2:=I_{2n+1}+E_{n+1,n+1}(1)\in\Mat_{2n+1}(k)$. Then, for $j=1,2$, we have
\begin{align*}
g_j^T Qg_j&=\begin{pmatrix} I_{2n+1}&0\\P_j^T&I_{2n+1}\end{pmatrix}\cdot \begin{pmatrix} D_1&A\\0&D_2\end{pmatrix}\cdot\begin{pmatrix} I_{2n+1}&P_j\\0&I_{2n+1}\end{pmatrix}\\&=Q+\begin{pmatrix} 0&D_1 P_j\\P_j^T D_1&0\end{pmatrix}+\begin{pmatrix} 0&0\\0&P_j^T D_1 P_j+P_j^T A\end{pmatrix}.
\end{align*}
Note that $P_j^TA=A^TP_j$, and thus $P_j^T D_1 P_j+P_j^T A$ is symmetric. Then it follows from Lemma \ref{Lemma-XTAX_i,i=...} that for all $x=(x_1,\ldots,x_{4n+2})^T\in V_n$ we have
\[(g_jx)^TQ(g_jx)=x^TQx+\sum_{i=1}^{2n+1}(P_j^TD_1 P_j+P_j^T A)_{ii}\cdot x_{2n+1+i}^2.\]
So $(g_jx)^TQ(g_jx)=x^TQx$ if and only if $(P_j^T D_1 P_j+P_j^T A)_{ii}=0$, for all $i=1,\ldots,2n+1$. Note that $X_1^T=T_0(1)$ and $X_2^T=X_2$. 
As
\[(T_0(1)\cdot D_1\cdot T_2(1)+T_0(1)\cdot A)_{ii}=\begin{cases} 0& \text{if $i=1$}\\(D_1)_{i-1,i-1}+A_{i-1,i},& \text{if $i\in\{2,\ldots,2n+1\}$}\end{cases}\]
we get $(g_1x)^TQ(g_1x)=x^TQx$ if and only if $(D_1)_{i,i}=A_{i,i+1}$, for all $i=1,\ldots,2n$. Likewise we get $(g_2x)^TQ(g_2x)=x^TQx$ if and only if $(D_1)_{i,i}=A_{i,i}$, for all $i\in\{1,\ldots,2n+1\}{\setminus\{n+1\}}$. As $A_{i,i+1}=0=A_{i,i}$, for all $i=1,\ldots,n$, it follows that $Q$ is $G$-invariant if and only if
\begin{align*}
(1)&\ \text{$(D_1)_{i,i}=0$, for all $i=1,\ldots,n$, and }\\
(2)&\ \text{$(D_1)_{n+1,n+1}=A_{n+1,n+2}$ and}\\
(3)&\ \text{$(D_1)_{i,i}=A_{i,i+1}=A_{i,i}$, for all $i=n+2,\ldots,2n$ and}\\
(4)&\ \text{$(D_1)_{2n+1,2n+1}=A_{2n+1,2n+1}$.}
\end{align*}
All four conditions can be met if and only if $A_{i,i+1}=A_{i,i}$, for all $i=n+2,\ldots,2n$, that is, if and only if $\Omega\in\EuScript{F}'$. Also note that we have $D_1=D$, for the diagonal matrix $D\in\Mat_{2n+1}(k)$ as described in (\ref{EquationForD1}).\\
\hfill\\
\noindent {\bf{3. Isometry classes}}: Let $D'\in\Diag_{2n+1}(k)$ and let $Q_{\Omega}(D')$ be as defined in the statement. Also note that all assumptions of Lemma \ref{Lemma-M^TQ(D')M-Q(E_tt(eta_tt.x))}, which we prove in the appendix, are satisfied. If $\Omega=0$, then $D=0$. Thus $Q_{\Omega}(D')$ is isometric to $Q_{\Omega}(0)$, by Lemma \ref{Lemma-M^TQ(D')M-Q(E_tt(eta_tt.x))} (a)

Next let $\Omega\neq 0$. Then $D\neq 0$. Set $\eta:=R^{-1}DR^{-T}$. A quick matrix calculation reveals that $\eta=\widetilde{I}_{2n+1}D\widetilde{I}_{2n+1}$. Also there is a minimal $t\in\{1,\ldots,n+1\}$ such that $\eta_{t,t}\neq 0$. Now it follows from Lemma \ref{Lemma-M^TQ(D')M-Q(E_tt(eta_tt.x))} (b) that there is a unique $x\in\EuScript{P}$ such that $Q_{\Omega}(D')$ is isometric to $Q_{\Omega}(E_{t,t}(\eta_{t,t}^{-1}x))$. This completes the proof.\\
\qed\\\\

\subsection{$C_n(f)$}\label{subsectionCn(f)}
Let $f\in k[T]$ be an irreducible polynomial of degree $m$ over $k$ in the indeterminate $T$, let $n\geq 1$ be an integer and let $\Pi\in \Mat_{mn}(k)$ be the companion matrix of $f^n$. That means $\Pi$ has ones on the first sub-diagonal and the final column of $\Pi$ is $(a_0,ldots,a_{mn-1})^T$, where $f^n=a_0+a_1T+\ldots+a_{mn-1}T^{mn-1}+T^{mn}$. Then there is an indecomposable $kG$-module $C_n(f)$ given as follows. Let $C_n(f)$ be the $k$-vector space $k^{2mn}$ with the standard basis of unit coordinate column vectors $\EuScript{B}:=\{e_1,\ldots,e_{2mn}\}$. The $G$-action is given by 
\[g_1=\begin{pmatrix} I_{mn}&I_{mn}\\0&I_{mn}\end{pmatrix},\quad g_2=\begin{pmatrix} I_{mn}&\Pi\\0&I_{mn}\end{pmatrix}.\]
By $\op{C}_{\Mat_{mn}(k)}(\Pi)$ we denote the set of all elements in $\Mat_{mn}(k)$ that commute with $\Pi$. Also we define $\EuScript{W}:=\{B\in\Mat_{mn}(k):\ \Pi^TB=B \Pi\}$.
\begin{Lemma}\label{Lemma-End-kG(C_n(f))}
\[\End_{kG}(C_n(f))=\left\{\begin{pmatrix} A&B\\0&A\end{pmatrix}:\ A\in \op{C}_{\Mat_{mn}(k)}(\Pi), B\in\Mat_{mn}(k)\right\}.\]
\end{Lemma}

Proof: Let $M=\begin{pmatrix} A&B\\C&D\end{pmatrix}\in \Mat_{2mn}(k)$, where $A,B,C,D\in \Mat_{mn}(k)$. Then $M\in \End_{kG}(C_n(f))$ if and only if $Mg_1=g_1M$ and $Mg_2=g_2M$. The first equation holds if and only if $A=D$ and $C=0$. The second equation then holds if and only if $A\Pi=\Pi A$. This completes the proof.\\
\qed\\\\

\begin{Lemma}\label{Lemma-G-invariant forms on Cn(f)}
The $G$-invariant forms on $C_n(f)$ are precisely the matrices of the form
\[F=\begin{pmatrix} 0&B\\B&D\end{pmatrix},\]
where $B\in \EuScript{W}$ and $D\in\Mat_{mn}(k)$. Furthermore $F$ is symplectic if and only if $B$ is invertible and $D$ is symmetric and hollow.
\end{Lemma}

Proof: Let $F=\begin{pmatrix} A&B\\C&D\end{pmatrix}\in \Mat_{2mn}(k)$, where $A,B,C,D\in \Mat_{mn}(k)$. Then $F$ is a $G$-invariant form on $C_n(f)$ if and only if $g_1^TFg_1=F$ and $g_2^TFg_2=F$. One checks quickly that both equations are satisfied precisely if $F$ is as described in the statement.\\
\qed\\\\

For various reasons it is difficult to work within $\Mat_{2mn}(k)$ and instead it is more convenient to switch from $\Pi$ to its Jordan form. We set up this step in what follows. First let $\epsilon_1,\ldots,\epsilon_m$ be the roots of $f$ in a splitting field $E$. Note that $f$ has no repeated roots, as $k$ is perfect. Now $k[\epsilon_i]=\{x_0+x_1\epsilon_i+\ldots+x_{m-1}\epsilon_i^{m-1}: x_j\in k\}$ is a field extension over $k$ of degree $m$, for all $i\in\{1,\ldots,m\}$. Also note that $k[\epsilon_i]$ is perfect. Finally let $\sigma_i: k[\epsilon_1]\rightarrow k[\epsilon_i]$ be the unique $k$-homomorphism with $\sigma_i(\epsilon_1)=\epsilon_i$. In particular $\sigma_i$ is the identity map if restricted to $k$. For simplicity set $\epsilon:=\epsilon_1$.

Next consider the Jordan block matrix $J:=\diag(J_n(\sigma_1(\epsilon)),\ldots, J_n(\sigma_m(\epsilon)))\in \GL_{mn}(E)$, where $J_n(\alpha)=\alpha\cdot I_n+T_2(1)\in\GL_{mn}(E)$ denotes the Jordan matrix with $\alpha\in E$ on the main diagonal and ones on the super diagonal and zeros everywhere else.
\begin{Lemma}\label{Lemma-IfXJ=JX,thenX=diag(X1,...,XM)}
Let $X\in \Mat_{mn}(E)$. Then\\
\hfill\\
(1) $JX=XJ$ if and only if $X=\diag(X_1,\ldots,X_m)$, for upper triangular Toeplitz-matrices $X_1,\ldots,X_m\in\Mat_n(E)$.\\
\hfill\\
(2) $J^TX=XJ$ if and only if $X=\diag(X_1,\ldots,X_m)$, for lower triangular Hankel-matrices $X_1,\ldots,X_m\in\Mat_n(E)$. In particular $X$ is symmetric.
\end{Lemma}

Proof: Write $X=(X_{i,j})_{i,j=1,\ldots,m}$, where $X_{i,j}\in\Mat_{n}(E)$, for all $i,j=1,\ldots,m$.\\
\noindent (1) We have $JX=XJ$ if and only if $J_n(\epsilon_i) X_{i,j}=X_{i,j} J_n(\epsilon_j)$, for all $i,j=1,\ldots,m$. Let $Y\in\Mat_{n}(E)$ such that $J_n(\epsilon_j)Y=YJ_n(\epsilon_i)$. Then 
\[\epsilon_i\cdot Y+T_2(1)\cdot Y=J_n(\epsilon_i) Y=Y J_n(\epsilon_j)=\epsilon_j\cdot Y+Y\cdot T_2(1),\]
and so $(\epsilon_i+\epsilon_j)Y=YT_2(1)+T_2(1)Y$. By induction $(\epsilon_i+\epsilon_j)^tY=YT_2(1)^t+T_2(1)^tY$, for all integers $t\geq 1$. But $T_2(1)^n=0$, and thus $Y=0$ ensues, if $i\neq j$. If $i=j$, then $YT_2(1)=T_2(1)Y$. It is now an easy exercise to show that $Y$ is an upper triangular Toeplitz matrix.\\
\hfill\\
(2) Let $R:=\diag(\widetilde{I}_n,\ldots,\widetilde{I}_n)\in\Mat_{mn}(k)$. Then $J^T=R^{-1}JR$. Hence $J^TX=XJ$ if and only if $J(RX)=(RX)J$. Now the statement follows from part (1) since $\widetilde{I}_n Y$, is a lower triangular Hankel-matrix, for any upper triangular Toeplitz-matrix $Y$.\\
\qed\\\\

Note that $\Pi$ is similar to $J$. Thus there exists $V\in \GL_{mn}(E)$ so that $\Pi=V^{-1} J V$. If $V^{(i)}$ denotes the $i$-th column of $V$, then $V^{(i+1)}=V\cdot \Pi^{(i)}=J\cdot V^{(i)}$ and consequently $V^{(i+1)}=J^{i}\cdot V^{(1)}$, for all $i=1,\ldots, mn-1$. In particular $V$ is uniquely determined by its first column. 

Furthermore observe $(V^{(1)})_{tn}\neq 0$, for all $t\in\{1,\ldots,m\}$ as otherwise an inductive argument gives $(V^{(i+1)})_{tn}=(J\cdot V^{(i)})_{tn}=\epsilon_t\cdot (V^{(i)})_{tn}=0$, for all $i=1,\ldots,mn-1$, and thus the $tn$-th row of $V$ is zero in contradiction to $V$ being invertible. Now let us write $V^{(1)}=(v_1,\ldots,v_m)^T$, for row vectors $v_i\in E^n$. Then there are upper triangular Toeplitz matrices $T_{v_1},\ldots T_{v_n}$ such that $v_i=T_{v_i}\cdot e_n$, for all $i=1,\ldots,n$. Note that all $T_{v_i}$ are invertible and $\diag(T_{v_1},\ldots T_{v_n})$ commutes with $J$, by Lemma \ref{Lemma-IfXJ=JX,thenX=diag(X1,...,XM)}. Hence we can replace $V$ by $\diag(T_{v_1}^{-1},\ldots T_{v_n}^{-1})\cdot V$ and thus without lose of generality we may assume that $V^{(1)}=(\underbrace{e_n,\ldots,e_n}_{\text{$m$-times}})^T$.

Next let $r,s\in\{1,\ldots,mn\}$ and write $r=(t-1)n+q$, for $q\in\{1,\ldots,n\}$. Then, since $V^{(s)}=J^{s-1}\cdot V^{(1)}$, we have
\[V_{r,s}=\left(J^{s-1}\cdot V\right)_{r,1}=\left[(J_n(\sigma_t(\epsilon)))^{s-1}\right]_{q,n}=\sigma_t\left((J_n(\epsilon)^{s-1})_{q,n}\right),\]
and so $V_{r,s}=\sigma_t(V_{q,s})$. In particular there is some $V'\in \Mat_{n\times (mn)}(k[\epsilon])$ such that
\begin{equation}\label{euqationOfV}
V=\begin{pmatrix} V'\\\sigma_2(V')\\\vdots\\ \sigma_m(V')\end{pmatrix}.
\end{equation}
Throughout the remainder of the paper we set for any $X\in\Mat_n(k[\epsilon])$
\begin{equation}\label{equation-DX,VX}
\begin{split}
\EuScript{D}(X)&:=\diag(X,\sigma_2(X)\ldots,\sigma_m(X))\\
\EuScript{V}(X)&:=V^T\cdot \diag(X,\sigma_2(X)\ldots,\sigma_m(X))\cdot V
\end{split}
\end{equation}
\begin{Lemma}\label{LemmaV^{-T}AV^{-1}andVC_Mat_{mn}(k)(Pi)V^{-1}}
(1) $\EuScript{W}=\{\EuScript{V}(X): X\in \Mat_{n}(k[\epsilon])\text{ a lower triangular Hankel-matrix}\}$\\
\hfill\\
(2) $\op{C}_{\Mat_{mn}(k)}(\Pi)=\{V^{-1}\EuScript{D}(X)V: X\in \Mat_{n}(k[\epsilon])\text{ an upper triangular Toeplitz-matrix}\}$.
\end{Lemma}

Proof: (1) "$\Leftarrow$": Let $A=\EuScript{V}(X)$, for some lower triangular Hankel-matrix $X\in \Mat_n(k[\epsilon])$. Then $J^T(V^{-T}AV^{-1})=(V^{-T}AV^{-1})J$, by Lemma \ref{Lemma-IfXJ=JX,thenX=diag(X1,...,XM)}(2) and thus $\Pi^TA=A\Pi$. It remains to show that $A\in\Mat_{mn}(k)$. In (\ref{euqationOfV}) we write $V'=(V_1\ldots V_m)\in\Mat_{n\times (mn)}(k[\epsilon])$, where $V_t\in\Mat_n(k[\epsilon])$. Also write $A=(A_{i,j})_{i,j=1,\ldots,m}$, 
for $A_{i,j}\in\Mat_n(E)$. Then 
\[A_{i,j}=\sum_{t=1}^m \sigma_t(V_i^TXV_j).\]
Now for every $k$-automorphism $\tau$ on $E$ we have $\tau(A_{i,j})=\sum_{t=1}^m (\tau\circ\sigma_t)(V_i^TXV_j)$. But $\tau$ permutes the set $\{\sigma_1,\ldots,\sigma_m\}$ and thus $\tau(A_{i,j})=A_{i,j}$. Hence $A_{i,j}\in\Mat_n(k)$ and thus $A\in\Mat_{mn}(k)$. In particular $A\in\EuScript{W}$.

"$\Rightarrow$": First let $R:=\diag(\widetilde{I}_n,\ldots,\widetilde{I}_n)\in\Mat_{mn}(k)$. Then $B:=V^TRV\in\EuScript{W}$, by the previous paragraph. Note that $B$ is invertible and thus $V^{-1}=B^{-1}V^TR$. Next let $A\in \EuScript{W}$. Then $J^T(V^{-T}AV^{-1})=(V^{-T}AV^{-1})J$. Hence by Lemma \ref{Lemma-IfXJ=JX,thenX=diag(X1,...,XM)}(2) we have $V^{-T}AV^{-1}=\diag(X_1,\ldots,X_m)$, for lower triangular Hankel-matrices $X_1,\ldots,X_m\in\Mat_{n}(E)$. Hence $(R^TV)\cdot (B^{-T}AB^{-1})\cdot (V^TR)=\diag(X_1,\ldots,X_n)$. Now with (\ref{euqationOfV}) it follows that $X_t=\sigma_t\left((\widetilde{I}_nV')\cdot (B^{-T}AB^{-1}) \cdot (V'^T\widetilde{I}_n)\right)$, for all $t=1,\ldots,m$. As $(\widetilde{I}_nV')\cdot (B^{-T}AB^{-1}) \cdot (V'^T\widetilde{I}_n)\in\Mat_n(k[\epsilon])$ we have $X_1\in\Mat_n(k[\epsilon])$ and $X_t=\sigma_t(X_1)$, for all $t=1,\ldots,m$.\\
\hfill\\
(2) Still let $B\hspace*{-0.1cm}=\hspace*{-0.1cm}V^TRV\hspace*{-0.1cm}\in\EuScript{W}$. Then $A\in\op{C}_{\Mat_{mn}(k)}(\Pi)$ if and only if $BA\hspace*{-0.1cm}\in\EuScript{W}$. By part (1) this is equivalent to $BA=\EuScript{V}(X)$, for some lower triangular Hankel-matrix $X\in \Mat_{n}(k[\epsilon])$, that is, $A=V^{-1}\diag(\widetilde{I}_nX,\sigma_2(\widetilde{I}_nX)\ldots,\sigma_m(\widetilde{I}_nX))V$, since $B^{-1}V^T=V^{-1}R^{-1}$ and $R^{-1}=R$. As $\widetilde{I}_nX$ is an upper triangular Toeplitz-matrix if and only if $X$ is a lower triangular Hankel-matrix, the statement follows.\\
\qed\\\\

\begin{Lemma}\label{Lemma-AB=B^TA and BC=CB}
Let $A\in\EuScript{W}$ and let $B,C\in\op{C}_{\Mat_{mn}(k)}(\Pi)$. Then $AB=B^TA$ and $BC=CB$.
\end{Lemma}

Proof: By Lemma \ref{LemmaV^{-T}AV^{-1}andVC_Mat_{mn}(k)(Pi)V^{-1}}(1) elements in $\EuScript{W}$ are symmetric and one checks quickly that $AB\in\EuScript{W}$. Hence $AB=(AB)^T=B^TA^T=B^TA$. Moreover $BC=CB$ follows from Lemma \ref{LemmaV^{-T}AV^{-1}andVC_Mat_{mn}(k)(Pi)V^{-1}}(2) since upper triangular Toeplitz matrices commute.\\
\qed\\\\

\begin{Lemma}\label{Lemma-End(Cn(f))/J(End(Cn(f)))}
We have $\End_{kG}(C_n(f))/J(\End_{kG}(C_n(f)))\cong k[\epsilon]$.
\end{Lemma}

Proof: By Lemmas \ref{Lemma-End-kG(C_n(f))} and \ref{LemmaV^{-T}AV^{-1}andVC_Mat_{mn}(k)(Pi)V^{-1}} the elements in $\End_{kG}(C_n(f))$ are precisely the matrices $M=\begin{pmatrix} A&B\\0&A\end{pmatrix}$, where $A=V^{-1}\EuScript{D}(X)V$, for an upper triangular Toeplitz-matrix $X\in \Mat_{n}(k[\epsilon]$ and $B\in\Mat_{mn}(k)$. Note that $M$ is invertible if $X_{1,1}$ is non-zero, and $M$ is nilpotent otherwise. Now one checks that $\End_{kG}(C_n(f))\rightarrow k[\epsilon]: M\mapsto X_{1,1}$ is an epimorphism of $k$-algebras with kernel $J(\End_{kG}(C_n(f)))$.\\
\qed\\\\

The following two lemmas are of a technical nature and we use them in our main theorems below.
\begin{Lemma}\label{Lemma-leftcosetsOfHinUT_n(K)}
Let $K$ be a perfect field of characteristic two and $n\geq 1$ an integer. For $r\in\{1,\ldots,n\}$ we define
\begin{equation}\label{equation-HrK}
\EuScript{H}_r(K)=\left\{H_{n+r}(1_K)+H_{n+r+1}(1_K) A^2: A\in\EuScript{UT}_n(K)\right\}\subseteq \EuScript{LH}_n(K).
\end{equation}
Next let $A\in\EuScript{LH}_n(K)$ such that $\rk(A)=n+1-r$. Then for every integer $s\geq 0$ there is a unique $C\in\EuScript{H}_{r+2s}(K)$ such that $AB^2=C$, for some $B\in\EuScript{UT}_n(K)$. Furthermore if $s=0$, then $B$ is invertible.
\end{Lemma}

Proof: By assumption $A=\sum_{i=r}^n H_{n+i}(x_i)$, where $x_i\in K$, for all $i=r,\ldots,n$. There are $y_i,z_i\in K$ such that $y_i^2=x_{r+2i-2}$, for all $i$ such that $1\leq i\leq \frac{n-r+2}{2}$ and $y_i=0$ otherwise, and $z_i^2=x_{r+2i-1}$, for all $i$ such that $1\leq i\leq \frac{n-r+1}{2}$ and $z_i=0$ otherwise. Now set $A_1=\sum_{i=1}^n T_i(y_i)$ and $A_2=\sum_{i=1}^n T_i(z_i)$. Then $A_1,A_2\in\EuScript{UT}_n(K)$ with $A_1$ invertible and $A=H_{n+r}(1_K)\cdot A_1^2+H_{n+r+1}(1_K)\cdot A_2^2$. 

Next for $B\in\EuScript{UT}_n(K)$ we have $AB^2=H_{n+r}(1_K)\cdot (A_1B)^2+H_{n+r+1}(1_K)\cdot (A_2B)^2$. Note that if $B=A_1^{-1}T_{s+1}(1_K)$, then $AB^2\in \EuScript{H}_{r+2s}(K)$. More general $AB^2\in \EuScript{H}_{r+2s}(K)$ if and only if $H_{n+r}(1_K)\cdot (A_1B)^2=H_{n+r+2s}(1_K)$. But in this case $H_{n+r}(1_K)\cdot B^2=H_{n+r+2s}(1_K)\cdot A_1^{-2}$ and thus $H_{n+r+1}(1_K)\cdot B^2=H_{n+r+1}(1_K)A_1^{-2}$. As then $AB^2=H_{n+r+2s}(1_K)+H_{n+r+2s+1}(1_K)\cdot (A_1^{-1}A_2)^2$ uniqueness ensues. 

Finally if $B$ is not invertible, then $\rk(AB^2)<\rk(A)$. Consequently if $s=0$, then $B$ must be invertible.\\
\qed\\\\

\begin{Lemma}\label{Lemma-Pi^TchiPi_i,i=Pi^Tchi_i,i-ThenThreeCases}
Let $\D \chi=H_{n+r}(1_k)+\sum_{\substack{i=1\\\text{$i$ is odd}}}^{n-r} H_{n+r+i}(\chi_i)\in\Mat_n(k[\epsilon])$, for some $r\in\{1,\ldots,n+1\}$. Here $H_j(\star)\in\Mat_n(k[\epsilon])$ is defined as the $n\times n$ -- zero matrix if $j>2n$. Then $(\Pi^T \EuScript{V}(\chi) \Pi)_{i,i}=(\Pi^T\EuScript{V}(\chi))_{i,i}$, for all $i=1,\ldots,mn$, if and only if $\chi$ satisfies one of the following cases:
\begin{enumerate}
\item $\chi=0$
\item $n+r$ is even, $f\in\{T,T+1_k\}$
and $\chi=\begin{cases} H_{2n}(1_k)&\text{if $r=n$}\\H_{n+r}(1_k)+H_{n+r+1}(1_k)&\text{if $r<n$}\end{cases}$
\item $n+r$ is odd, $f\not\in\{T,T+1_k\}$ 
and $\chi_{i}=(\epsilon+\epsilon^2)^{-\frac{i+1}{2}}$, for all odd $i\in\{1,\ldots,n-r\}$.
\end{enumerate}
\end{Lemma}

Proof: We claim that $(\Pi^T \EuScript{V}(\chi) \Pi)_{i,i}=(\Pi^T\EuScript{V}(\chi))_{i,i}$, for all $i=1,\ldots,mn$, if and only if
\begin{equation}\label{Equation-JalphaTchiJalpha_ii}
(J(\epsilon)^T\chi J(\epsilon))_{i,i}=(J(\epsilon)^T\chi)_{i,i}\text{, for all $i=1,\ldots,n$.}
\end{equation}
As $\EuScript{V}(\chi)\in\EuScript{W}$ we have $\Pi^T\EuScript{V}(\chi)(\Pi+I_{mn})\in\EuScript{W}$ is symmetric. Hence $(\Pi^T \EuScript{V}(\chi) \Pi)_{i,i}=(\Pi^T\EuScript{V}(\chi))_{i,i}$, for all $i=1,\ldots,mn$ if and only if $\Pi^T\EuScript{V}(\chi)(\Pi+I_{mn})$ is hollow, which by Lemma \ref{Lemma-XTAX_i,i=...} holds if and only if $V^{-T}\Pi^T\EuScript{V}(\chi)(\Pi+I_{mn})V^{-1}$ is hollow. Since $\Pi=V^{-1}JV$, we have $V^{-T}\Pi^T\EuScript{V}(\chi)(\Pi+I_{mn})V^{-1}=J^TV^{-T}\EuScript{V}(\chi) V^{-1}(J+I_{mn})=J^T\cdot \EuScript{D}(\chi)\cdot (J+I_{mn})$. This matrix is hollow if and only if $J(\epsilon)^T\chi(J(\epsilon)+I_n)$ is hollow. So the claim follows.

Clearly (\ref{Equation-JalphaTchiJalpha_ii}) holds if $\chi=0$. Hence let now $\chi\neq 0$, that is, $r\leq n$. Then
\[(J(\epsilon)^T\chi J(\epsilon))_{i,i}=\sum_{s=1}^n \chi_{s,s}\cdot (J(\epsilon)_{s,i})^2,\]
by Lemma \ref{Lemma-XTAX_i,i=...}, and 
\[(J(\epsilon)^T\chi)_{i,i}=\sum_{s=1}^n\chi_{s,i}\cdot J(\epsilon)_{s,i}.\]
Hence (\ref{Equation-JalphaTchiJalpha_ii}) holds if and only if $\chi_{i,i}(\epsilon+\epsilon^2)=\chi_{i-1,i-1}+\chi_{i-1,i}$, for all $i=1,\ldots, n$, where $\chi_{0,0}$ and $\chi_{0,1}$ are defined as zero. Since $\chi_{i,j}=0$, whenever $i+j<n+r$, it follows that (\ref{Equation-JalphaTchiJalpha_ii}) is equivalent to
\begin{equation}\label{Equation-chi_11-and-chi_ii}
\chi_{i,i}(\epsilon+\epsilon^2)=\chi_{i-1,i-1}+\chi_{i-1,i}\ \text{, for all integers $i$ with $\frac{n+r}{2}\leq i\leq n$.}
\end{equation}
First assume that $n+r$ is even. Then
\[\chi_{i,i}=\begin{cases} 1_k,&\text{if $i=\frac{n+r}{2}$}\\0_k,&\text{otherwise.}\end{cases}\]
For $j=\frac{n+r}{2}$ we get from (\ref{Equation-chi_11-and-chi_ii}) that
\[(\epsilon+\epsilon^2)=\chi_{j,j}(\epsilon+\epsilon^2)=\chi_{j-1,j-1}+\chi_{j-1,j}=0.\]
Hence $\epsilon\in\{0_k,1_k\}$, that is, $f\in\{T,T+1_k\}$. Also (\ref{Equation-chi_11-and-chi_ii}) is now equivalent to $\chi_{i-1,i-1}=\chi_{i-1,i}$, for all $i=\frac{n+r+2}{2},\ldots,n$. This condition is empty if $r=n$, and thus $\chi=H_{2n}(1_k)$ is the only solution. So let $r<n$. Note that as $\chi$ is a Hankel-matrix we have for all $i,j$ such that $i+j\geq n+r$ that $\chi_{i,j}=\chi_{i+j-n,n}=\chi_{r+(i+j-(n+r)),n}=\chi_{i+j-(n+r)}$. Thus (\ref{Equation-chi_11-and-chi_ii}) is equivalent to $\chi_{2(i-1)-(n+r)}=\chi_{2i-1-(n+r)}$, for all $i\geq \frac{n+r+2}{2}$, that is, $\chi_{s}=\chi_{s+1}$, for all even $s\in\{0,\ldots,n-r-2\}$. Since by assumption $\chi_0=1_k$ and $\chi_s=0$, for all even $s\in\{2,\ldots,n-r-2\}$, we conclude that in this case (\ref{Equation-chi_11-and-chi_ii}) holds if and only if $\chi=H_{n+r}(1_k)+H_{n+r+1}(1_k)$.

Next assume that $n+r$ is odd. Then
\[\chi_{i,i}=\begin{cases} 0_k,&\text{if $i<\frac{n+r}{2}$}\\\chi_{2i-(n+r)},&\text{if $i>\frac{n+r}{2}$.}\end{cases}\]
Then (\ref{Equation-chi_11-and-chi_ii}) is equivalent to $\chi_1(\epsilon+\epsilon^2)=1_k$ and $\chi_s(\epsilon+\epsilon^2)=\chi_{s-2}$, for all odd $s\in\{3,\ldots,n-r\}$. In particular in this case (\ref{Equation-chi_11-and-chi_ii}) holds if and only if $\epsilon\not\in\{0_k,1_k\}$ (that is, $f\not\in\{T,T+1_k\}$) and $\chi_s=(\epsilon+\epsilon^2)^{-\frac{s+1}{2}}$, for all odd $s\in\{1,\ldots,n-r\}$.\\
\qed\\\\

Next let $\EuScript{H}=\D\left\{H_{n+1}(1_k)+\hspace*{-0.2cm}\sum_{\substack{i=2\\\text{$i$ is even}}}^n\hspace*{-0.2cm} H_{n+i}(a_i)\in \Mat_n(k[\epsilon]): a_i\in k[\epsilon]\right\}$. Note that $\EuScript{H}=\EuScript{H}_1(k[\epsilon])$, as defined in (\ref{equation-HrK}). Also for every $A\in\EuScript{H}$ let $S(A):=\begin{pmatrix} 0&\EuScript{V}(A)\\\EuScript{V}(A)&0\end{pmatrix}$. We can now state and prove the main Theorem of this section.
\begin{Theorem}\label{Theorem-C_n(f)}
(a) {\bf{symplectic forms:}} We have $\EuScript{IS}(C_n(f))=\{S(A):\ A\in\EuScript{H}\}$. In particular if $k$ is finite, then $|\EuScript{IS}(C_n(f))|=|k|^{m\cdot \lfloor\frac{n}{2}\rfloor}$.\\
\hfill\\
\noindent (b) {\bf{quadratic forms:}} Let $A\in\EuScript{H}$. There is a $G$-invariant quadratic form with associated bilinear form $S(A)$, precisely in the following cases:
\begin{enumerate}
\item $n$ is odd, $f\in\{T,T+1_k\}$ and $A=\begin{cases} 1_k,&\text{if $n=1$}\\H_{n+1}(1_k)+H_{n+2}(1_k),& \text{if $n\geq 3$,}\end{cases}$
\item $n\geq 2$ is even, $f\not\in\{T,T+1_k\}$ and $a_i=(\epsilon+\epsilon^2)^{-\frac{i}{2}}$, for all even $i\in\{2,\ldots,n\}$.
\end{enumerate}
In this case $\{Q_A(D): D\in\Diag_{mn}(k)\}$ gives all such quadratic forms, where $Q_A(D):=\begin{pmatrix} D_A&\EuScript{V}(A)\\0&D\end{pmatrix}$ and $D_A\in\Mat_{mn}(k)$ is the diagonal matrix with $(D_A)_{i,i}=(\EuScript{V}(A))_{i,i}$, for all $i=1,\ldots,mn$. Finally $\EuScript{Q}(S(A))=\{Q_{A}(E_{t,t}(v^{-1} x)):\ x\in\EuScript{P}\}$, where $v:=((\EuScript{V}(A))^{-1})_{t,t}\neq 0$, for some $t\in\{1,\ldots,mn\}$.
\end{Theorem}

Proof: (a) {\bf{1. Notation}}: By Lemma \ref{Lemma-G-invariant forms on Cn(f)} the indecomposable symplectic $kG$-modules $(C_n(f),S)$ are precisely given by $S(B,D):=\begin{pmatrix} 0&B\\B&D\end{pmatrix}$, where $B\in \EuScript{W}$ is invertible and $D\in\Mat_{mn}(k)$ is symmetric and hollow. Let $S=S(B,D)$ be such a form. Also let $M:=\begin{pmatrix} X&Y\\0&X\end{pmatrix}\in\Aut_{kG}(C_n(f))$, that is, $X\in\op{C}_{\GL_{mn}(k)}(\Pi)$ and $Y\in\Mat_{mn}(k)$. Then
\[M^TSM=\begin{pmatrix} 0&X^TBX\\X^TBX&Y^TBX+X^TBY+X^TDX\end{pmatrix}.\]

\noindent {\bf{2. $S$ is isometric to $S(A)$, for a unique $A\in\EuScript{H}$}}: First note that there is some $Z\in\Mat_{mn}(k)$ such that $D=Z+Z^T$. Hence with $X=I_{mn}$ and $Y=B^{-1}Z$ we get $M^TSM=S(B,0)$, that is, $S$ is isometric to $S(B,0)$. Likewise we see that $S$ is isometric to some $S(C,0)$, if and only if $C=X^TBX$, for some $X\in \op{C}_{\GL_{mn}(k)}(\Pi)$. 

Next note that $X^TBX=BX^2$, by Lemma \ref{Lemma-AB=B^TA and BC=CB}. Furthermore by Lemma \ref{LemmaV^{-T}AV^{-1}andVC_Mat_{mn}(k)(Pi)V^{-1}} we have $X=V^{-1}\EuScript{D}(X')V$, for some invertible $X'\in\EuScript{UT}_n(k[\epsilon])$ and $B=\EuScript{V}(B')$, for some invertible $B'\in\EuScript{LH}_n(k[\epsilon])$. By Lemma \ref{Lemma-leftcosetsOfHinUT_n(K)} there is a unique $A\in\EuScript{H}_1(k[\epsilon])=\EuScript{H}$ such that $A=B'(X')^2$, for some invertible $X'\in\EuScript{UT}_n(k[\epsilon])$. Now $BX^2=\EuScript{V}(B'(X')^2)=\EuScript{V}(A)$, that is, $S$ is isometric to $S(A)$. \\

(b) {\bf{1. $G$-invariant quadratic forms}}: Let $A\in\EuScript{H}$. Then the quadratic forms with associated bilinear form $S(A)$ are precisely the matrices $\begin{pmatrix} D_1&\EuScript{V}(A)\\0&D_2\end{pmatrix}$, where $D_1,D_2\in\Diag_{mn}(k)$. Let $Q$ be such a matrix. Note that $Q$ is $G$-invariant if and only if $x^T(g_j^TQg_j+Q)x=0$, for all $x\in C_n(f)$ and $j=1,2$. One checks quickly that for $j=1,2$
\[g_j^TQg_j+Q=\begin{pmatrix} 0&D_1P_j\\ P_j^T D_1&P_j^TD_1P_j+P_j^T\EuScript{V}(A)\end{pmatrix},\]
where $P_1=I_{mn}$ and $P_2=\Pi$. Clearly $D_1$, $D_1+\EuScript{V}(A)$ and $\Pi^TD_1\Pi$ are all symmetric. Also $\Pi^T\EuScript{V}(A)$ is symmetric, since $V\Pi=JV$ and $AJ=J^TA$. Therefore $g_j^TQg_j+Q$ is symmetric and so $Q$ is $G$-invariant if and only if $P_j^TD_1P_j+P_j^T$ is hollow for $j=1,2$, that is, $(D_1)_{i,i}=\EuScript{V}(A)_{i,i}$ and $(\Pi^TD_1\Pi)_{i,i}=(\Pi^T\EuScript{V}(A))_{i,i}$, for all $i=1,\ldots,mn$. In particular note that if such $D_1$ exists, it is uniquely determined by the main diagonal of $\EuScript{V}(A)$. Furthermore $D_2$ has no bearing on the $G$-invariance of $Q$.

We need to inquire when such a $D_1$ exists. Assume that $(D_1)_{i,i}=\EuScript{V}(A)_{i,i}$, for all $i=1,\ldots,mn$. Then $D_1+\EuScript{V}(A)$ is symmetric and hollow and thus by Lemma \ref{Lemma-XTAX_i,i=...} we have $(\Pi^TD_1\Pi)_{i,i}=(\Pi^T\EuScript{V}(A)\Pi)_{i,i}$, for all $i=1,\ldots,mn$. Hence $Q$ is $G$-invariant if and only if $(\Pi^T\EuScript{V}(A)\Pi)_{i,i}=(\Pi^T\EuScript{V}(A))_{i,i}$, for all $i=1,\ldots,mn$. Applying Lemma \ref{Lemma-Pi^TchiPi_i,i=Pi^Tchi_i,i-ThenThreeCases} gives all instances of $A$ such that $Q$ is $G$-invariant.\\
\hfill\\
\noindent {\bf{2. isometry classes}}: Let $A$ be such that $Q_A(D)$ is a $G$-invariant quadratic form with associated bilinear form $S(A)$ and set $\eta:=\EuScript{V}(A)^{-1} D_A \EuScript{V}(A)^{-T}$. First note that $A$ is not hollow. This is true because if $n$ is odd, then $A_{\frac{n+1}{2},\frac{n+1}{2}}=1_k$ and if $n$ is even, then $A_{\frac{n+2}{2},\frac{n+2}{2}}=(\epsilon+\epsilon^2)^{-1}$. So it follows from Lemma \ref{Lemma-XTAX_i,i=...} that $D_A$ and consequently $\eta$ are not hollow. Thus there is some $t\in\{1,\ldots, mn\}$ such that $\eta_{t,t}\neq 0$. Since $(D_A)_{i,i}=\EuScript{V}(A)_{i,i}$, for all $i\in\{1,\ldots,mn\}$ it follows from Lemma \ref{Lemma-XTAX_i,i=...} that $\eta_{t,t}=(\EuScript{V}(A)^{-1} \EuScript{V}(A) \EuScript{V}(A)^{-T})_{t,t}=(\EuScript{V}(A)^{-1})_{t,t}$. Hence for $v:=\eta_{t,t}$ Lemma \ref{Lemma-M^TQ(D')M-Q(E_tt(eta_tt.x))}(b) shows that $Q_A(D)$ is isometric to $Q_A(E_{t,t}(v^{-1}x))$, for a unique $x\in\EuScript{P}$.\\
\qed\\\\

\subsection{$C_n(f)\oplus C_n(f)^*$}\label{subsection-Cn(f)^2}
We use the same notation as in subsection \ref{subsectionCn(f)}. In particular $f=f(T)$, $\Pi$, $\EuScript{W}$, $\epsilon$ and $V$ are all the same as before. As recall the notations $\EuScript{V}(-)$ and $\EuScript{D}(-)$ as given in (\ref{equation-DX,VX}). First we show that $C_n(f)$ is self-dual.
\begin{Lemma}\label{Lemma-C_n(f)*=C_n(f)}
$C_n(f)^*\cong C_n(f)$
\end{Lemma}

Proof: As every matrix is similar to its transpose, there is some $A\in\GL_{mn}(k)$ such that $\Pi^T=A^{-1}\Pi A$. Hence with $B=\begin{pmatrix} A&0\\ 0&A\end{pmatrix}$ we get $B^{-1}g_1^TB=g_1$ and $B^{-1}g_2^TB=g_2$. In particular $C_n(f)$ is self-dual as claimed.\\
\qed\\\\

For $\varphi,\mu,\psi\in\Mat_n(k[\epsilon])$ set $S(\varphi,\psi,\mu):=\begin{pmatrix} 0&\sigma\\ \sigma&0\end{pmatrix}$, where $\sigma=\begin{pmatrix} \EuScript{V}(\varphi)&\EuScript{V}(\psi)\\ \EuScript{V}(\psi)&\EuScript{V}(\mu)\end{pmatrix}$. Furthermore let $\EuScript{K}$ be the set of all triplets $(\varphi,\psi,\mu)$ such that: $\varphi,\mu,\psi$ are lower triangular Hankel-matrices $\Mat_n(k[\epsilon])$ such that $\mu=H_{n+r+2s+1}(1_k)$ and
\[\varphi=H_{n+r}(1_k)+\sum_{\substack{i=1\\ \text{$i$ odd}}}^{2s-1}H_{n+r+i}(\varphi_i),\ \psi=H_{n+1}(1_k)+\sum_{i=1}^tH_{n+i}(\psi_i),\]
for integers $r,s$ such that $2\leq r\leq n+1$ and $0\leq s\leq \frac{n-r+1}{2}$, and where $t:=\left\lfloor\frac{n-r-2s+1}{2}\right\rfloor$ and $\psi_1\neq 1_k$, if $t\geq 1$. Here $H_j(-)$ is defined as the zero $n\times n$--matrix if $j>2n$.

In the following we work with respect to the basis \[\EuScript{B}:=\{e_1\oplus 0,\mydots,e_{mn}\oplus 0,0\oplus e_1,\mydots,0\oplus e_{mn},e_{mn+1}\oplus 0,\mydots,e_{2mn}\oplus 0,0\oplus e_{mn+1},\mydots,0\oplus e_{2mn}\}.\]

\begin{Theorem}\label{Theorem-C_n(f)+C_n(f)*}
(a) {\bf{symplectic forms:}} We have $\EuScript{IS}(C_n(f)^2)=\{S(\varphi,\mu,\psi): (\varphi,\mu,\psi)\in\EuScript{K}\}$. Furthermore if $k$ is finite, then
\[|\EuScript{IS}(C_n(f)^2)|=\begin{cases} n|k|^{m\left(\frac{n-2}{2}\right)},&\text{if $n$ is even}\\ \frac{n+1}{2}\cdot |k|^{m\left(\frac{n-1}{2}\right)}+\frac{n-1}{2}\cdot |k|^{m\left(\frac{n-3}{2}\right)},&\text{if $n$ is odd.}\end{cases}\]
\hfill\\
\noindent (b) {\bf{quadratic forms:}} Let $(\varphi,\psi,\mu)\in \EuScript{K}$. There is a $G$-invariant quadratic form with associated bilinear form $S=S(\varphi,\psi,\mu)$ precisely if $\mu=0$ and one of the following cases hold:
\begin{enumerate}
\item $\varphi=0$
\item $n+r$ is even, $f\in\{T,T+1\}$, $\varphi_1=1_k$, if $s\geq 1$ and $\varphi_i=0_k$, for all odd $i\in\{3,\ldots,2s-1\}$
\item $n+r$ is odd, $f\not\in\{T,T+1\}$ and $\varphi_{i}=(\epsilon+\epsilon^2)^{-\frac{i+1}{2}}$, for all odd $i\in\{1,\ldots,2s-1\}$.
\end{enumerate}
In this case $\{Q_{\sigma}(D'): D'\in\Diag_{2mn}(k)\}$ gives all such quadratic forms, where $Q_{\sigma}(D'):=\begin{pmatrix} D & \sigma\\ 0& D'\end{pmatrix}\in\Mat_{4mn}(k)$ and $D\in\Mat_{2mn}(k)$ is the diagonal matrix such that $D_{i,i}=\sigma_{i,i}$, for all $i=1,\ldots,2mn$. 

If $\varphi=\mu=0$, then $\EuScript{Q}(S(\varphi,\psi,\mu))=\{Q_{\sigma}(0)\}$. Otherwise there is $t\in\{1,\ldots,2mn\}$ so that $\eta_{t,t}\neq 0$, where $\eta:=\sigma^{-1}D\sigma^{-T}$ and $\EuScript{Q}(S(\varphi,\psi,\mu))=\{Q_{\sigma}(E_{t,t}(\eta_{t,t}^{-1}x)):x\in\EuScript{P}\}$.
\end{Theorem}

Proof: (a) {\bf{1. Notation}}: By Lemma \ref{Lemma-End(Cn(f))/J(End(Cn(f)))} Remark \ref{Remark-symplectic_forms_on_M^2} applies. Hence with respect to the basis $\{e_1\oplus 0,\mydots,e_{2mn}\oplus 0,0\oplus e_1,\mydots,0\oplus e_{2mn}\}$ the indecomposable symplectic $kG$-modules $(C_n(f)^2,S)$ are precisely given by $S=\begin{pmatrix} S_1&S_2\\S_2^T&S_3\end{pmatrix}$, where all $S_i$ are $G$-invariant forms on $C_n(f)$, while $S_1$ and $S_3$ are also alternating and degenerate. By Lemmas \ref{Lemma-G-invariant forms on Cn(f)} and \ref{LemmaV^{-T}AV^{-1}andVC_Mat_{mn}(k)(Pi)V^{-1}} we have $S_i=\begin{pmatrix}0&\EuScript{V}(\omega_i)\\\EuScript{V}(\omega_i)&\Omega_i\end{pmatrix}$, where all $\omega_i\in\Mat_n(k[\epsilon])$ are lower triangular Hankel-matrices and $\Omega_i\in\Mat_{mn}(k)$ such that $\Omega_1$ and $\Omega_3$ are symmetric and hollow, and $\omega_1$ and $\omega_3$ are singular. In the following let $S$ be such a form.\\
\hfill\\
\noindent {\bf{2. Change of basis}}: From now on we work with respect to the basis $\EuScript{B}$. Then $S=\begin{pmatrix}0&\omega\\\omega&\Omega\end{pmatrix}$, where $\omega:=\begin{pmatrix} \EuScript{V}(\omega_1)&\EuScript{V}(\omega_2)\\\EuScript{V}(\omega_2)&\EuScript{V}(\omega_3)\end{pmatrix}$ and $\Omega:=\begin{pmatrix} \Omega_1&\Omega_2\\\Omega_2^T&\Omega_3\end{pmatrix}$. Observe that $\omega$ is self-transpose and since $S$ is non-degenerate $\omega$ is invertible. Consequently $\omega_2$ must be invertible. Also note that $(C_n(f)^2,S(\varphi,\psi,\mu))$ is an indecomposable symplectic $kG$-module for any triplet $(\varphi,\psi,\mu)$ of lower triangular Hankel-matrix in $\Mat_n(k[\epsilon])$, where $\varphi$ and $\mu$ and singular and $\psi$ is invertible.\\
\hfill\\
\noindent {\bf{3. $S$ is isometric to $S(\omega_1,\omega_2,\omega_3)$}}: Next let $M\in\Aut_{kG}(C_n(f)^2)$. Then by Lemma \ref{Lemma-End-kG(C_n(f))} and with respect to basis $\EuScript{B}$ we have $M=\begin{pmatrix}\upsilon&\Upsilon\\0&\upsilon\end{pmatrix}$, where $\upsilon:=\begin{pmatrix} a&b\\c&d\end{pmatrix}$, for $a,b,c,d\in \op{C}_{\Mat_{mn}(k)}(\Pi)$, and $\Upsilon\in\Mat_{2mn}(k)$. Furthermore
\[M^TSM=\begin{pmatrix}\upsilon^T&0\\\Upsilon^T&\upsilon^T\end{pmatrix}\cdot\begin{pmatrix}0&\omega\\\omega&\Omega\end{pmatrix}\cdot \begin{pmatrix}\upsilon&\Upsilon\\0&\upsilon\end{pmatrix}=\begin{pmatrix}0&\upsilon^T\omega\upsilon\\ \upsilon^T\omega\upsilon&\Upsilon^T\omega\upsilon+\upsilon^T\omega \Upsilon+\upsilon^T\Omega\upsilon\end{pmatrix}.\]
As $\Omega$ is symmetric and hollow there is some $Z\in\Mat_{2mn}(k)$ such that $\Omega=Z+Z^T$. So if we set $\upsilon=I_{2mn}$ and $\Upsilon=\omega^{-1}Z$, then $M$ is invertible and $M^TSM=S(\omega_1,\omega_2,\omega_3)$. In particular $S$ is isometric to $S(\omega_1,\omega_2,\omega_3)$. Also note that $S$ is isometric to some $S(\varphi,\psi,\mu)$ if and only if $\upsilon^T\omega\upsilon=\begin{pmatrix} \EuScript{V}(\varphi)&\EuScript{V}(\psi)\\ \EuScript{V}(\psi)&\EuScript{V}(\mu)\end{pmatrix}$, for some $\upsilon$ has described above.\\
\hfill\\
\noindent {\bf{4. The isometry class of $S(\omega_1,\omega_2,\omega_3)$}}: We need to understand $\upsilon^T\omega\upsilon$. First set $\omega_i':=\EuScript{V}(\omega_i)$, for $i=1,2,3$. Also by Lemma \ref{LemmaV^{-T}AV^{-1}andVC_Mat_{mn}(k)(Pi)V^{-1}} we have $a=V^{-1}\EuScript{D}(\alpha) V$, $b=V^{-1}\EuScript{D}(\beta) V$, $c=V^{-1}\EuScript{D}(\gamma) V$ and $d=V^{-1}\EuScript{D}(\delta)V$, for upper triangular Toeplitz-matrices $\alpha,\beta,\gamma,\delta\in\Mat_n(k[\epsilon])$. (Note that $M$ is invertible if and only if $\upsilon$ is invertible if and only if $\alpha\delta+\beta\gamma$ is invertible.) Then
\begin{align*}
\upsilon^T\omega\upsilon&=\begin{pmatrix} a^T&c^T\\b^T&d^T\end{pmatrix}\cdot \begin{pmatrix} \omega_1'&\omega_2'\\\omega_2'&\omega_3'\end{pmatrix}\cdot \begin{pmatrix} a&b\\c&d\end{pmatrix}\\
&=\begin{pmatrix} \omega_1'a^2+\omega_3'c^2&\omega_1'ab+\omega_2'bc+\omega_2'ad+\omega_3'cd\\\omega_1'ab+\omega_2'bc+\omega_2'ad+cd&\omega_1'b^2+\omega_3'd^2\end{pmatrix}\\
&=\begin{pmatrix} \EuScript{V}(\omega_1\alpha^2+\omega_3 \gamma^2)&\EuScript{V}(\omega_1 \alpha\beta+\omega_2\alpha\delta+\omega_2\gamma\beta+\omega_3\gamma\delta)\\\EuScript{V}(\omega_1 \alpha\beta+\omega_2\alpha\delta+\omega_2\gamma\beta+\omega_3\gamma\delta)&\EuScript{V}(\omega_1 \beta^2+\omega_3\delta^2)\end{pmatrix},
\end{align*}
where the second equality uses Lemma \ref{Lemma-AB=B^TA and BC=CB}. Consequently the isometry class of $S(\omega_1,\omega_2,\omega_3)$ contains precisely all $S(\varphi,\psi,\mu)$ where
\begin{equation}\label{EquationsLittleOmega}
\begin{split}
\varphi&=\omega_1\alpha^2+\omega_3 \gamma^2\\
\mu&=\omega_1 \beta^2+\omega_3\delta^2\\
\psi&=\omega_1 \alpha\beta+\omega_2\alpha\delta+\omega_2\gamma\beta+\omega_3\gamma\delta
\end{split}
\end{equation}
for some upper triangular Toeplitz-matrices $\alpha,\beta,\gamma,\delta\in\Mat_n(k[\epsilon])$ such that $\alpha\delta+\beta\gamma$ is invertible. \\
\hfill\\
\noindent {\bf{5. Every isometry class contains some $S(\varphi,\psi,\mu)$ with $(\varphi,\psi,\mu)\in\EuScript{K}$}}: In the following let $\EuScript{C}$ denote a given isometry class. First let $\EuScript{C}$ contains $S(0,\psi,0)$, for some invertible $\psi\in\EuScript{LH}_n(k[\epsilon])$. Note that $\psi^{-1}\widetilde{I}_n\in \EuScript{UT}_n(k[\epsilon])$ and set $\alpha=\psi^{-1}\widetilde{I}_n$, $\beta=\gamma=0$ and $\delta=I_n$. Then $S(0,\psi,0)$ is isometric to $S(0,\widetilde{I}_n,0)$. Note that $(0,\widetilde{I}_n,0)\in\EuScript{K}$, with $r=n+1$, $s=0$ and $t=0$.

Next assume that $\EuScript{C}$ contains some $S(\varphi,\psi,\mu)$ such that one of $\varphi$ and $\mu$ equals zero but not both. Note that $S(\varphi,\psi,\mu)$ is isometric to $S(\mu,\psi,\varphi)$ with $\alpha=\delta=I_n$ and $\beta=\gamma=0$. Hence we may assume that $\mu=0$. By Lemma \ref{Lemma-leftcosetsOfHinUT_n(K)} there is some $A\in\EuScript{UT}_n(k[\epsilon])$ such that $\varphi A^2\in\EuScript{H}_r(k[\epsilon])$, for some integer $r\in\{2,\ldots,n\}$. Also note that $\psi^{-1}\widetilde{I}_n\in\EuScript{UT}_n(k[\epsilon])$. Hence with $\alpha=A$, $\beta=\gamma=0$ and $\delta=\alpha^{-1}\psi^{-1}\widetilde{I}_n$ we obtain that $S(\varphi,\psi,\mu)$ is isometric to $S(\varphi A^2,\widetilde{I}_n,0)$. Note that $(\varphi A^2,\widetilde{I}_n,0)\in\EuScript{K}$ with $r\in\{2,\ldots,n\}$ and $s=\lfloor \frac{n-r+1}{2}\rfloor$.

Finally we can assume that if $S(\varphi,\psi,\mu)\in\EuScript{C}$, then $\varphi,\mu\neq 0$. Without lose of generality we may choose $S(\varphi,\psi,\mu)\in\EuScript{C}$ such that for $r:=n+1-\rk(\varphi)$ and $j:=n+1-\rk(\mu)$ we have $j\geq r$ and $j-r$ is maximal in $\EuScript{C}$. Note that $r,j\in\{2,\ldots,n\}$. If $j-r=2i$, for some integer $i\geq 0$, then set $\alpha=I_n$, $\beta=T_{i+1}(x)$, where $x\in k[\epsilon]$ such that $x^2=(\varphi_{r,n})^{-1}\cdot\mu_{j,n}$, $\gamma=0$ and $\delta=I_n$. Now $S(\varphi,\psi,\mu)$ is isometric to $S(\varphi,\varphi\beta+\psi,\varphi\beta^2+\mu)$. A quick calculation reveals that $\varphi\beta^2+\mu$ has a rank of at most $\rk(\mu)-1$ thus contradicting the maximality of $j-r$ in $\EuScript{C}$. In particular we have $j=r+2s+1$, for some integer $s$ such that $0\leq s\leq \lfloor\frac{n-r-1}{2}\rfloor$. 

Recall $\EuScript{H}_r(k[\epsilon])\subseteq\EuScript{LH}_n(k[\epsilon])$ as defined in (\ref{equation-HrK}). By Lemma \ref{Lemma-leftcosetsOfHinUT_n(K)} there are invertible $A,D\in\EuScript{UT}_n(k[\epsilon])$ such that $\varphi A^2\in\EuScript{H}_r(k[\epsilon])$ and $\mu D^2\in\EuScript{H}_j(k[\epsilon])$. Hence with $\alpha=A,\beta=\gamma=0$ and $\delta=D$ we get that $S(\varphi,\psi,\mu)$ is isometric to $S(\varphi A^2,\psi AD,\mu D^2)$. In particular we may assume that $\varphi\in\EuScript{H}_r(k[\epsilon])$ and $\mu\in\EuScript{H}_j(k[\epsilon])$. Also since $j+1-r$ is even it follows from Lemma \ref{Lemma-leftcosetsOfHinUT_n(K)} that there is some $B\in\EuScript{UT}_n(k[\epsilon])$ such that $\varphi B^2\in\EuScript{H}_{j+1}(k[\epsilon])$. Now $\varphi B^2=H_{n+j+1}(1_k)+H_{n+j+2}(1_k)\cdot X_1^2$ and $\mu=H_{n+j}(1_k)+H_{n+j+1}(1_k)\cdot X_3^2$, for some $X_1,X_3\in \EuScript{UT}_n(k[\epsilon])$. Note that $\varphi (BX_3)^2+H_{n+j}(1_k)\cdot (I_n+T_2(1_k)X_1)^2=\mu$. Hence set $\alpha=I_n$, $\beta=BX_3$, $\gamma=0$ and $\delta=I_n+T_2(1_k)X_1$. Then $\alpha\delta+\beta\gamma=I_n+T_2(1_k)X_1$ is invertible and $S(\varphi,(\varphi\beta+\psi)\delta^{-1},H_{n+j}(1_k))$ is isometric to $S(\varphi,\psi,\mu)$. In particular we may assume that $\mu=H_{n+j}(1_k)$.

Since $\varphi\in\EuScript{H}_r(k[\epsilon])$ we have $\varphi=\varphi'+\varphi''$, where
\[\varphi'=H_{n+r}(1_k)+\sum_{\substack{i=1\\\text{$i$ is odd}}}^{2s-1} H_{n+r+i}(\varphi_{r+i})\quad\text{ and }\quad \varphi''=\sum_{\substack{i=2s+1\\\text{$i$ is odd}}}^{n-r} H_{n+r+i}(\varphi_{r+i}),\]
with $\varphi_i\in k[\epsilon]$, for all $i=r+1,\ldots,n$. Note that $\varphi''=H_{n+j}(1_k)\cdot C^2$, where $C=\sum_{i=1}^n T_i(c_i)$, with $c_i^2=\varphi_{r+2s+2i-1}$, for all $i=1,\ldots,\lfloor\frac{n+1-r-2s}{2}\rfloor$ and $c_i=0$ otherwise. Now with $\alpha=\delta=I_n$, $\beta=0$ and $\gamma=C$ we get that $S(\varphi,\psi,H_{n+j}(1_k))$ is isometric to $S(\varphi+H_{n+j}(1_k) C^2,\psi+H_{n+j}(1_k) C,H_{n+j}(1_k))$. Since $\varphi+H_{n+j}(1_k) C^2=\varphi'$ we may assume that $\varphi=\varphi'$, as required in the statement. 

Next set $t:=\lfloor\frac{n-r-2s+1}{2}\rfloor$ and write $\psi=\psi'+\psi''$, where $\psi'=\sum_{i=1}^t H_{n+i}(\psi_i)$ and $\psi''=\sum_{i=t+1}^n H_{n+i}(\psi_i)$, for $\psi_i\in k[\epsilon]$, for all $i=1,\ldots,n$. Note that $t\geq 1$, as $r+2s+1\leq n$, and $\psi_1\neq 0$, as $\psi$ is invertible. Set $\alpha=I_n$, $\beta=\gamma=0$ and $\delta=I_n+\psi^{-1}\psi''$. Then $S(\varphi,\psi,H_{n+j}(1_k))$ is isometric to $S(\varphi,\psi\delta,H_{n+j}(1_k)\delta^2)=S(\varphi,\psi+\psi'',H_{n+j}(1_k)+H_{n+j}(1_k)(\psi^{-1}\psi'')^2)$. Note that $\psi+\psi''=\psi'$ as required. Furthermore note that $\psi^{-1}\psi''$ is an upper triangular Toeplitz-matrix of rank equal to $\rk(\psi'')=n-t$. Consequently $\rk(H_{n+j}(1_k)(\psi^{-1}\psi'')^2)=n+1-(2t+j)$, by (\ref{equation-rank(AB)}). But $t=\lfloor \frac{n-j+2}{2}\rfloor$. Hence $2t\geq n-j+1$ and so $n+1-(2t+j)\geq 0$. Thus $H_{n+j}(1_k)(\psi^{-1}\psi'')^2=0$. Overall we have found some $S(\varphi,\psi,\mu)\in\EuScript{C}$ where $(\varphi,\psi,\mu)\in\EuScript{K}$ with $r\in\{2,\ldots,n-1\}$, $0\leq s\leq \frac{n-r+1}{2}$ and $t=\lfloor\frac{n-r-2s+1}{2}\rfloor$.\\
\hfill\\
\noindent {\bf{6. Uniqueness}}: Let $(\varphi,\psi,\mu), (\varphi',\psi',\mu')\in\EuScript{K}$ with accompanying $r,r'\in\{2,\ldots,n+1\}$, and integers $s,s'$ with $0\leq s\leq \frac{n-r+1}{2}$ and $0\leq s'\leq \frac{n-r'+1}{2}$ so that $S(\varphi,\psi,\mu)$ and $S(\varphi',\psi',\mu')$ are isometric. Then there are upper triangular Toeplitz-matrices $\alpha=\sum_{i=1}^n T_i(\alpha_i)$, $\beta=\sum_{i=1}^n T_i(\beta_i)$, $\gamma=\sum_{i=1}^n T_i(\gamma_i)$ and $\delta=\sum_{i=1}^n T_i(\delta_i)$ such that $\varphi'=\varphi\alpha^2+\mu \gamma^2$, $\mu'=\varphi\beta^2+\mu\delta^2$ and $\psi'=\varphi\alpha\beta+\psi\alpha\delta+\psi\gamma\beta+\mu\gamma\delta$ and $\alpha\delta+\beta\gamma$ is invertible. Note that the latter holds precisely if $\alpha_1\delta_1+\beta_1\gamma_1\neq 0$. 

First suppose that $r\neq r'$. Without lose of generality $r'<r$, that is, $\rk(\varphi)<\rk(\varphi')$. Then
\[\rk(\varphi)<\rk(\varphi')=\rk(\varphi\alpha^2+\mu\gamma^2)\leq \min\{\rk(\varphi\alpha^2),\rk(\mu\gamma^2)\}\leq \rk(\varphi).\]
This contradiction forces $r=r'$. Next suppose that $s\neq s'$ and without lose of generality $s'<s$. Then $\rk(\mu)<\rk(\mu')$. In particular $\mu'\neq 0$. Also since $\rk(\mu\delta^2)\leq\rk(\mu)<\rk(\mu')$ it follows that $\mu'$ and $\mu'+\mu\delta^2$ have the same rank. Finally since $\mu'+\mu\delta^2=\varphi\beta^2$ we get
\begin{align*}
n+1-(r+2s'+1)&=\rk(\mu')=\rk(\varphi\beta^2)\overset{(\ref{equation-rank(AB)})}{=}\rk(\varphi)+2\rk(\beta)-2n\\&=n+1-(r+2n-2\rk(\beta)).
\end{align*}
As now $2s'+1=2n-2\rk(\beta)$ we have a contradiction and thus $s=s'$. In particular it follows that $\mu=\mu'$ and $t=t'$.

If $r=n+1$, then $\varphi=0=\varphi'$, $\mu=0=\mu'$ and $\psi=I_n=\psi'$. So we may assume that $r\leq n$ and thus $\varphi,\varphi'\neq 0$. Next note that $\varphi\alpha^2=\varphi'+\mu\gamma^2$. Since $\rk(\varphi'+\mu\gamma^2)=r$ it follows that $\alpha$ has full rank and thus is invertible. Also observe that $\varphi\alpha^2=\varphi'+\mu\gamma^2\in\EuScript{H}_r(k[\epsilon])$. Since $\varphi=\varphi\cdot (I_n)^2 \in\EuScript{H}_r(k[\epsilon])$ the uniqueness statement in Lemma \ref{Lemma-leftcosetsOfHinUT_n(K)} implies that $\varphi=\varphi\alpha^2$. Consequently $\varphi=\varphi'+\mu\gamma^2$ and therefore $\varphi=\varphi'$ and $\mu\gamma^2=0$. The latter implies that $\rk(\gamma)\leq n-\frac{\rk(\mu)}{2}=\frac{n+r+2s-2}{2}$, and so $\gamma_i=0$, for all $i\in\{1,\ldots,\lfloor\frac{n-r+1-2s}{2}\rfloor\}$.

Also $0=\varphi\cdot (I_n+\alpha^2)=\varphi\cdot (I_n+\alpha)^2$ implies that $I_n+\alpha$ is not invertible, and so $\alpha_1=1_k$. Moreover $\rk(I_n+\alpha)\leq n-\frac{\rk(\varphi)}{2}=\frac{n+r-1}{2}$, and thus $\alpha_i=0$, for all $i\in\{2,\ldots\lfloor\frac{n-r+2}{2}\rfloor\}$. Note that $\varphi\beta^2=\mu (I_n+\delta)^2$. If $\varphi\beta^2\neq 0$, then
\begin{align*}
\rk(\varphi\beta^2)&=\rk(\varphi)+2\rk(\beta)-n=2\rk(\beta)+1-r\\
\rk(\mu (I_n+\delta)^2)&=\rk(\mu)+2\rk(I_n+\delta)-n=2\rk(I_n+\delta)-2s-r
\end{align*}
As this is a contradiction we must have $\varphi\beta^2=0$ and so $\rk(\beta)\leq n-\frac{\rk(\varphi)}{2}=\frac{n+r-1}{2}$. Hence $b_i=0$, for all $i\in \{1,\ldots,\lfloor\frac{n-r+2}{2}\rfloor\}$. Also $\mu (I_n+\delta)^2=0$, now implies that $I_n+\delta$ is not invertible, and so $\delta_1=1_k$. Moreover $\rk(I_n+\delta)\leq n-\frac{\rk(\mu)}{2}=\frac{n+r+2s-2}{2}$. Hence $\delta_i=0$, for all $i\in\{2,\ldots,\lfloor \frac{n-r-2s+1}{2}\rfloor\}$.

In summary we have $\alpha_1=\delta_1=1_k$, $\beta_1=\gamma_1=0$ and $\alpha_i=\beta_i=\gamma_i=\delta_i=0$, for all $i\in\{2,\ldots,t\}$. Next let $l\in\{1,\ldots,t\}$. Then $(\varphi\alpha\beta)_{l,n}=(\mu\gamma\delta)_{l,n}=(\psi\beta\gamma)_{l,n}=0$. Hence
\[\psi'_{l,n}=(\psi\alpha\delta)_{l,n}=\sum_{i=n+1-l}^n\sum_{j=i}^n \psi_{l,i}\cdot \alpha_{j-i+1}\cdot \delta_{n-j+1}=\psi_{l,n}.\]
Therefore $\psi'=\psi$ and thus overall we have established that $(\varphi,\psi,\mu)=(\varphi',\psi',\mu')$.\\
\hfill\\
\noindent {\bf{7. number of isometry classes if $k$ is finite}}: Let $(\varphi,\psi,\mu)\in\EuScript{K}$ and let $r,s,t,\varphi_i$ and $\psi_i$ be accordingly. For a fixed $r\in\{2,\ldots,n+1\}$ set $r':=n+1-r$ and $l:=\left\lfloor\frac{r'}{2}\right\rfloor$. Then $s\in \{0,1,\ldots,l\}$ and $t=l-s$. If $s=l$, that is, $t=0$, then for each odd $i\in\{1,\ldots,2s-1\}$ we have $|k|^m$ choices for $\varphi_i\in k[\epsilon]$. Hence there are $|k|^{ml}$ different classes. If $s<l$, that is, $t\geq 1$, then we have $|k|^{ms}$ different choices for $\varphi$, $|k|^m-1$ different choices for $\psi_1$ and $|k|^m$ different choices for all $\psi_2,\ldots,\psi_t$. Thus there are $(|k|^m-1)\cdot |k|^{m(s+t-1)}=(|k|^m-1)\cdot |k|^{m(l-1)}=|k|^{ml}-|k|^{m(l-1)}$ different classes. Overall for the number $x_{r'}$ of different isometry classes for a fixed $r$ we have
\begin{align*}
x_{r'}&=|k|^{ml}+\sum_{s=0}^{l-1} |k|^{ml}-|k|^{m(l-1)}=|k|^{ml}+l\cdot (|k|^{ml}-|k|^{m(l-1)})\\&=(l+1)|k|^{ml}-l|k|^{m(l-1)}=|k|^{-m}\cdot ((l+1)|k|^{m(l+1)}-l|k|^{ml})
\end{align*}
Note that $x_{r'}=x_{r'+1}$ for even $r'$. Moreover the positive part of $x_{r'}$ cancels with the negative part of $x_{r'+2}$. Also $x_0$ has no negative part. As $x=\sum_{r'=0}^{n-1}x_{r'}$ we see that $x$ is equal to the positive part of $x_{n-2}$ and $x_{n-1}$ survive. Finally if $n$ is even we have
\[x=2\cdot |k|^{-m}\cdot \left(\frac{n-2}{2}+1\right)|k|^{m\left(\frac{n-2}{2}+1\right)}=n|k|^{m\left(\frac{n-2}{2}\right)}.\] 
If $n$ is odd, then 
\begin{align*}
x&=|k|^{-m}\cdot \left(\left(\frac{n-3}{2}+1\right)|k|^{m\left(\frac{n-3}{2}+1\right)}+\left(\frac{n-1}{2}+1\right)|k|^{m\left(\frac{n-1}{2}+1\right)} \right)\\&=\left(\frac{n-1}{2}\right)|k|^{m\left(\frac{n-3}{2}\right)}+\left(\frac{n+1}{2}\right)|k|^{m\left(\frac{n-1}{2}\right)}. 
\end{align*}

(b) {\bf{1. $G$-invariant quadratic forms}}: Let $S=S(\varphi,\psi,\mu)$, where $(\varphi,\psi,\mu)\in\EuScript{K}$. Then the quadratic forms with associated bilinear from $S$ are precisely the matrices $\widehat{S}+D$, where $D\in\Diag_{4mn}(k)$. Let $Q$ be such a form and write $D=\diag(D_1,D_2)$, for $D_1,D_2\in\Diag_{2mn}(k)$. Then $Q=\begin{pmatrix} D_1&\sigma\\0&D_2\end{pmatrix}$, where $\sigma=\begin{pmatrix} \EuScript{V}(\varphi)&\EuScript{V}(\psi)\\ \EuScript{V}(\psi)&\EuScript{V}(\mu)\end{pmatrix}$. Also we have $g_j=\begin{pmatrix} I_{2mn}&P_j\\0&I_{2mn}\end{pmatrix}$, for $j=1,2$, where $P_1=I_{2mn}$ and $P_2=\begin{pmatrix} \Pi &0\\0&\Pi\end{pmatrix}$. Note that $Q$ is $G$-invariant if and only if $x^T(g_j^TQg_j+Q)x=0$, for all $x\in C_n(f)^2$ and $j=1,2$. One checks that 
\[g_j^TQg_j+Q=\begin{pmatrix} 0&D_1P_j\\P_j^TD_1&P_j^TD_1P_j+P_j^T\sigma\end{pmatrix}\]
Note that $(P_j^T\sigma)^T=P_j^T\sigma$. This is trivial if $j=1$ as $\sigma^T=\sigma$, and it is true for $j=2$, since $\EuScript{V}(\varphi),\EuScript{V}(\mu),\EuScript{V}(\psi)\in\EuScript{W}$, by Lemma \ref{LemmaV^{-T}AV^{-1}andVC_Mat_{mn}(k)(Pi)V^{-1}}(1). In particular $g_j^TQg_j+Q$ is symmetric for $j=1,2$. Hence $Q$ is $G$-invariant if and only if $P_j^TD_1P_j+P_j^T\sigma$ is hollow for $j=1,2$, that is, $(D_1)_{i,i}=\sigma_{i,i}$ and $(P_2^TD_1P_2)_{i,i}=(P_2^T\sigma)_{i,i}$, for all $i=1,\ldots,2mn$. Note that if such $D_1$ exists, it is uniquely determined by the main diagonal of $\sigma$ and there is no restriction on $D_2$. Next write $D_1=\diag(d_1,d_2)$, for diagonal matrices $d_1,d_2\in\Mat_{mn}(k)$, such that $(D_1)_{i,i}=\sigma_{i,i}$, for all $i=1,\ldots,2mn$, that is, $(d_1)_{i,i}\cong (\EuScript{V}(\varphi))_{i,i}$ and $(d_2)_{i,i}\cong (\EuScript{V}(\mu))_{i,i}$, for all $i=1,\ldots,mn$. Then $Q$ is $G$-invariant if and only if for all $i=1,\ldots,mn$ 
\begin{align*}
(\Pi^Td_1\Pi)_{i,i}=(\Pi^T\EuScript{V}(\varphi))_{i,i}\ \text{ and }\ (\Pi^Td_2\Pi)_{i,i}=(\Pi^T\EuScript{V}(\mu))_{i,i}.
\end{align*}
As $\EuScript{V}(\varphi),\EuScript{V}(\mu)\in\EuScript{W}$ are symmetric it follows from Lemma \ref{Lemma-XTAX_i,i=...} that $(\Pi^Td_1\Pi)_{i,i}=(\Pi^T\EuScript{V}(\varphi)\Pi)_{i,i}$ and $(\Pi^Td_2\Pi)_{i,i}=(\Pi^T\EuScript{V}(\mu)\Pi)_{i,i}$, for all $i=1,\ldots,mn$. In particular $Q$ is $G$-invariant if and only if $(\Pi^T\EuScript{V}(\varphi)\Pi)_{i,i}=(\Pi^T\EuScript{V}(\varphi))_{i,i}$ and $(\Pi^T\EuScript{V}(\mu)\Pi)_{i,i}=(\Pi^T\EuScript{V}(\mu))_{i,i}$, for all $i=1,\ldots,mn$. Note that Lemma \ref{Lemma-Pi^TchiPi_i,i=Pi^Tchi_i,i-ThenThreeCases} discusses which $\varphi$ and $\mu$ satisfy this condition. If $\mu=0$, then we get the three possibilities for $\varphi$ as stated in the theorem. If $\mu\neq 0$, then $r+2s+1\leq n$, and so $\varphi\neq 0$. Since $r+2s+1$ and $r$ have different parity we see quickly that none of the remaining possibilities for $\mu$ and $\varphi$ as given in Lemma \ref{Lemma-Pi^TchiPi_i,i=Pi^Tchi_i,i-ThenThreeCases} are feasible. In particular $Q$ is $G$-invariant precisely if $\mu$ and $\varphi$ are as stated and in this case $Q\cong Q_{\sigma}(D')$, for some diagonal $D'\in\Mat_{2mn}(k)$.\\ 
\hfill\\
\noindent {\bf{2. isometry classes}}: Next let $D'\in\Mat_{2mn}(k)$ be a diagonal matrix. If $\varphi=0$ and $\mu=0$, then $D=0$. Then $Q_{\sigma}(D')$ is isometric to $Q_{\sigma}(0)$, by Lemma \ref{Lemma-M^TQ(D')M-Q(E_tt(eta_tt.x))} (a). Hence let $\varphi\neq 0$ or $\mu\neq 0$. Then $\varphi$ and $\mu$ are not both hollow. So $\EuScript{V}(\varphi)$ and $\EuScript{V}(\mu)$ are not both hollow by Lemma \ref{Lemma-XTAX_i,i=...}. Thus $\sigma$ is not hollow and consequently $D\neq 0$. Set $\eta:=\sigma^{-1}D\sigma^{-T}$ and choose $t\in\{1,\ldots,2mn\}$ such that $\eta_{t,t}\neq 0$. Then by Lemma \ref{Lemma-M^TQ(D')M-Q(E_tt(eta_tt.x))} (b) there is a unique $x\in\EuScript{P}$ such that $Q_{\sigma}(D')$ is isometric to $Q_{\sigma}(E_{t,t}(\eta_{t,t}^{-1}x))$. In particular this completes the proof.\\
\qed\\\\

\subsection{$C_n(\infty)$ and $C_n(\infty)\oplus C_n(\infty)^*$}

Let $n\geq 1$ be an integer. Then $C_n(\infty)=k^{2n}$ and with respect to the basis $\EuScript{B}:=\{e_1,\ldots,e_{2n}\}$ the $G$-action on $C_n(\infty)$ is given by 
\[g_1=\begin{pmatrix} I_n&T_2(1)\\0&I_n\end{pmatrix}\ \text{and}\ g_2=\begin{pmatrix} I_n&I_n\\0&I_n\end{pmatrix}.\]
Note that if we switch to the basis $\EuScript{B}':=\{e_n,\ldots,e_1,e_{2n},\ldots,e_{n+1}\}$, then $g_1$ and $g_2$ act on $C_n(\infty)$ as $g_2$ and $g_1$, respectively, act on $C_n(f)$, where $f(T)=T\in k[T]$. Hence $C_n(\infty)$ and $C_n(T)$ have the same symplectic and quadratic forms and the same isometry classes. Note that with $f(T)=T$ we have $m=1$, $\epsilon=0_k$ and $V=\widetilde{I}_n$. Finally note that $\begin{pmatrix} \widetilde{I}_n&0\\0&\widetilde{I}_n\end{pmatrix}$ switches between the bases $\EuScript{B}$ and $\EuScript{B}'$. 

Also let $\EuScript{H}=\D\left\{H_{n+1}(1_k)+\hspace*{-0.2cm}\sum_{\substack{i=2\\\text{$i$ is even}}}^n\hspace*{-0.2cm} H_{n+i}(a_i)\in \Mat_n(k): a_i\in k\right\}$, and for every $A\in\EuScript{H}$ let $S(A):=\begin{pmatrix} 0&A\\A&0\end{pmatrix}$. With respect to the basis $\EuScript{B}$ we get by Theorem \ref{Theorem-C_n(f)}

\begin{Theorem}\label{Theorem-C_n(infty)}
(a) {\bf{symplectic forms}}: We have $\EuScript{IS}(C_n(\infty))=\{S(A): A\in \EuScript{H}\}$. In particular if $k$ is finite, then $|\EuScript{IS}(C_n(\infty))|=|k|^{\lfloor\frac{n}{2}\rfloor}$.\\
\hfill\\
\noindent (b) {\bf{quadratic forms}}: Let $A\in\EuScript{H}$. There is a $G$-invariant quadratic form with associated bilinear form $S(A)$ precisely if $n$ is odd and $A=1_k$, if $n=1$ and $A=H_{n+1}(1_k)+H_{n+2}(1_k)$ if $n\geq 3$. In this case $\{Q_A(D): D\in\Diag_n(k)\}$ gives all such quadratic forms, where $Q_A(D):=\begin{pmatrix} D_A&A\\0&D\end{pmatrix}$, and $D_A\in\Mat_{n}(k)$ is the diagonal matrix such that $(D_A)_{i,i}=A_{i,i}$, for all $i=1,\ldots,n$. Finally $\EuScript{Q}(S(A))=\{Q_{A}(E_{t,t}(v^{-1} x)):\ x\in\EuScript{P}\}$, where $v:=(A^{-1})_{t,t}\neq 0$, for some $t\in\{1,\ldots,n\}$.
\end{Theorem}

Since $C_n(T)$ is self-dual by Lemma \ref{Lemma-C_n(f)*=C_n(f)}, so is $C_n(\infty)$. For $\varphi,\mu,\psi\in\Mat_n(k)$ we define $S(\varphi,\psi,\mu):=\begin{pmatrix} 0&\sigma\\ \sigma&0\end{pmatrix}$, where $\sigma=\begin{pmatrix} \varphi&\psi\\ \psi&\mu\end{pmatrix}$. Also recall the definition of $\EuScript{K}$ given just before Theorem \ref{Theorem-C_n(f)+C_n(f)*}, (remember that now $k[\epsilon]=k$). Finally consider the basis $\EuScript{B}:= \{e_1\oplus 0,\ldots,e_n\oplus 0,0\oplus e_1,\ldots,0\oplus e_n,e_{n+1}\oplus 0,\ldots,e_{2n}\oplus 0,0\oplus e_{n+1},\ldots,0\oplus e_{2n}\}$. Then with Theorem \ref{Theorem-C_n(f)+C_n(f)*} we get
\begin{Theorem}\label{Theorem-C_n(infty)+C_n(infty)*}
(a) {\bf{symplectic forms}}: We have $\EuScript{IS}(C_n(f)^2)=\{S(\varphi,\mu,\psi): (\varphi,\mu,\psi)\in\EuScript{K}\}$. If $k$ is finite, then 
\[|\EuScript{IS}(C_n(f)^2)|=\begin{cases} n|k|^{\left(\frac{n-2}{2}\right)},&\text{if $n$ is even}\\ \frac{n+1}{2}\cdot |k|^{\left(\frac{n-1}{2}\right)}+\frac{n-1}{2}\cdot |k|^{\left(\frac{n-3}{2}\right)},&\text{if $n$ is odd.}\end{cases}\]
\noindent (b) {\bf{quadratic forms}}: Let $(\varphi,\mu,\psi)\in\EuScript{K}\}$. There is a $G$-invariant quadratic form with associated bilinear form $S(\varphi,\psi,\mu)$ precisely if $\mu=0$ and one of the following cases hold:
\begin{enumerate}
\item $\varphi=0$
\item $n+r$ is even, $\varphi_1=1_k$, if $s\geq 1$ and $\varphi_i=0_k$, for all odd $i\in\{3,\ldots,2s-1\}$
\end{enumerate}
In this case $\{Q_{\sigma}(D'): D'\in\Diag_{2n}(k)\}$ gives all such quadratic forms, where $Q_{\sigma}(D'):=\begin{pmatrix} D & \sigma\\ 0& D'\end{pmatrix}\in\Mat_{4n}(k)$, and $D\in\Mat_{2n}(k)$ is the diagonal matrix such that $D_{i,i}=\sigma_{i,i}$, for all $i=1,\ldots,2n$. 

If $\varphi=\mu=0$, then $\EuScript{Q}(S(\varphi,\psi,\mu))=\{Q_{\sigma}(0)\}$. Otherwise there is $t\in\{1,\ldots,2n\}$ so that $\eta_{t,t}\neq 0$, where $\eta:=\sigma^{-1}D\sigma^{-T}$, and $\EuScript{Q}(S(\varphi,\psi,\mu))=\{Q_{\sigma}(E_{t,t}(\eta_{t,t}^{-1}x)):\ x\in\EuScript{P}\}$.
\end{Theorem}

\section{Appendix}
We conclude the paper with the following result which was used in the proofs of Theorems \ref{Theorem_AnoplusBn}, \ref{Theorem-C_n(f)} and \ref{Theorem-C_n(f)+C_n(f)*} helping us enumerate the isometry classes of quadratic forms. Let $k$ be a perfect field of characteristic two and let $G$ be a group. Also let $V$ be a $kG$-module of dimension $2n$ and let $\EuScript{B}$ be a basis of $V$. Unless explicitly stated otherwise all matrices are with respect to the basis $\EuScript{B}$. Furthermore assume that
\[\left\{\begin{pmatrix} I_n&Y\\0&I_n\end{pmatrix}: y\in\Mat_n(k)\right\}\subseteq \Aut_{kG}(V).\]
Next let $(M,S)$ be symplectic where $S=\begin{pmatrix} 0&R\\R^T&0\end{pmatrix}$ for some $R\in\GL_n(k)$. For $F_1,F_2\in\Mat_n(k)$ set $Q_S(F_1,F_2):=\begin{pmatrix} F_1&R\\0&F_2\end{pmatrix}$. Note that for diagonal matrices $D_1,D_2\in\Mat_n(k)$ we have that $Q_S(D_1,D_2)$ corresponds to a quadratic form with associated bilinear form $S$. However in general $Q_S(D_1,D_2)$ may not be $G$-invariant. Finally let $\EuScript{P}$ be a full set of representatives for the distinct cosets of the additive subgroup $\{x^2+x:\ x\in k\}$ in $(k,+)$.

\begin{Lemma}\label{Lemma-M^TQ(D')M-Q(E_tt(eta_tt.x))}
With the above assumptions let $D_1,D_2\in\Mat_n(k)$ be diagonal matrices such that $Q_S(D_1,D_2)$ is $G$-invariant. 

(a) If $D_1=0$, then $Q_S(D_1,D_2)$ is isometric to $Q_S(0,0)$. 

(b) If $D_1\neq 0$ set $\eta:=R^{-1}D_1R^{-T}$ and choose $t\in\{1,\ldots,n\}$ such that $\eta_{t,t}\neq 0$. Then $Q_S(D_1,D_2)$ is isometric to $Q_S(D_1,E_{t,t}(\eta_{t,t}^{-1}x))$ for a  unique $x\in\EuScript{P}$.
\end{Lemma}

Proof: Let $M=\begin{pmatrix} I_n&Y\\0&I_n\end{pmatrix}\in\GL_{2n}(k)$, for $Y\in\Mat_n(k)$. Then by assumption $M\in\Aut_{kG}(V)$ and 
\begin{align*}
M^TQ_S(D_1,D_2)M&=\begin{pmatrix} D_1&D_1Y+R\\Y^TD_1&Y^TD_1Y+Y^TR+D_2\end{pmatrix}\\&\cong Q_S(D_1,Y^TD_1Y+Y^TR+D_2)
\end{align*}
First suppose that $D_1=0$. Then with $Y=R^{-T}D_2$ it follows that $M^TQ_S(D_1,D_2)M\cong Q_S(0,0)$. In particular part (a) of the statement holds.

Next assume that $D_1\neq 0$ and set $\eta:=R^{-1}D_1R^{-T}$. By Lemma \ref{Lemma-XTAX_i,i=...} there is some $t\in\{1,\ldots,n\}$ such that $\eta_{t,t}\neq 0$. In the following we construct a symmetric matrix $L\in\Mat_n(k)$ such that $E_{t,t}(\eta_{t,t}^{-1}x)+L+L\eta L\cong D_2$, for some $x\in\EuScript{P}$. First for all $i,j\in\{1,\ldots,n\}{\setminus\{t\}}$ we set $L_{i,j}=0$. Then by Lemma \ref{Lemma-XTAX_i,i=...} we have for all $i\in\{1,\ldots,n\}{\setminus\{t\}}$
\[(E_{t,t}(\eta_{t,t}^{-1}x)+L+L\eta L)_{i,i}=\eta_{t,t}\cdot (L_{t,i})^2.\]
Thus we can choose $L_{t,i}=L_{i,t}\in k$ so that $(E_{t,t}(\eta_{t,t}^{-1}x)+L+L\eta L)_{i,i}=(D_2)_{i,i}$. Also
\begin{align*}
((E_{t,t}(\eta_{t,t}^{-1}&x)+L+L\eta L)_{t,t}=\eta_{t,t}^{-1}x+L_{t,t}+\sum_{s=1}^{n} \eta_{s,s}\cdot (L_{s,t})^2\\&=\eta_{t,t}^{-1}\left(x+\eta_{t,t}\cdot L_{t,t}+\left(\eta_{t,t}\cdot (L_{t,t}\right)^2\right)+\sum_{\substack{s=1\\s\neq t}}^{n} \eta_{s,s}\cdot (L_{s,t})^2
\end{align*}
In particular there are $x\in\EuScript{P}$ and $L_{t,t}\in k$ so that $(E_{t,t}(\eta_{t,t}^{-1}x)+L+L\eta L)_{t,t}=(D_2)_{t,t}$. Thus 
$E_{t,t}(\eta_{t,t}^{-1}x)+L+L\eta L$ and $D_2$ have the same main diagonal and as $L$ is symmetric we get $D_2+L+L\eta L\cong E_{t,t}(\eta_{t,t}^{-1}x)$. Now set $Y=R^{-T}L$. Then $Y^TD_1Y+Y^TR+D_2=L\eta L+L+D_2\cong E_{t,t}(\eta_{t,t}^{-1}x)$. Consequently $M^TQ_S(D_1,D_2)M\cong Q_S(D_1,E_{t,t}(\eta_{t,t}^{-1}x))$. 

It remains to show that $x\in\EuScript{P}$ is unique. Let $q_x$ denote the $G$-invariant quadratic form that corresponds to $Q_S(D_1,E_{t,t}(\eta_{t,t}^{-1}x))$. First note that
\[\begin{pmatrix}R^{-1}&0\\0&I_n\end{pmatrix} \cdot Q_S(D_1,E_{t,t}(\eta_{t,t}^{-1}x)) \cdot \begin{pmatrix}R^{-T}&0\\0&I_n\end{pmatrix}=\begin{pmatrix}R^{-1}D_1R^{-T}&I_n\\0&E_{t,t}(\eta_{t,t}^{-1}x)\end{pmatrix}.\]
Hence there is a basis $\{u_1,\ldots,u_n,v_1,\ldots,v_n\}$ with respect to which $q_x$ corresponds to the matrix $\begin{pmatrix}R^{-1}D_1R^{-T}&I_n\\0&E_{t,t}(\eta_{t,t}^{-1}x)\end{pmatrix}$. Note that $\EuScript{B'}:=\{u_1,v_1,\ldots,u_n,v_n\}$ is a symplectic basis as with respect to $\EuScript{B'}$ the form $S$ corresponds to $\diag(\widetilde{I}_2,\ldots,\widetilde{I}_2)\in\GL_{2n}(k)$. Next set $\EuScript{A}_x:=\sum_{i=1}^n q_x(u_i)q_x(v_i)$. Then $\EuScript{A}_x=(R^{-1}D_1R^{-T})_{t,t}\cdot \eta_{t,t}^{-1}\cdot x=\eta_{t,t}\cdot \eta_{t,t}^{-1}\cdot x=x$. Finally assume that $q_x$ and $q_y$ are isometric for $x,y\in\EuScript{P}$. Then by [\cite{Grove}, Theorem 13.13] there is some $\delta\in k$ such that $\EuScript{A}_x=\EuScript{A}_y+\delta^2+\delta$. Hence $x=y+\delta^2+\delta$ and so $x=y$ ensues.\\
\qed\\\\


\begin{thebibliography}{99}

\bibitem{Conlon}
S. B. Conlon. Modular representations of $C_2\times C_2$. \emph{J. Austral. Math. Soc.}, {\bf{10}} (1969), 363-366.

\bibitem{GowWillems}
R. Gow, W. Willems. A Note on Green Correspondence and Forms. \emph{Com. in Algebra} {\bf{23}}(4) (1995), 1239-1248.

\bibitem{GowWillems2}
R. Gow, W. Willems. Quadratic Geometries, Projective Modules and Idempotents. \emph{J. Algebra} {\bf{160}} (1993), 257-272.

\bibitem{Grove}
L.C. Grove. \emph{Classical Groups and Geometric Algebra}, Graduate Studies in Mathematics. {\bf{39}},  AMS, (2001).

\bibitem{Wall}
G.E.Wall. On the conjugacy classes in the unitary, symplectic and orthogonal groups, \emph{J. Aus. Math. Soc.} {\bf{3}} (1963), 1-62.

\bibitem{WillemsThesis}
W. Willems. Metrische Modulen \"uber Gruppenringen. Phd-Thesis, Johannes Gutenberg-Universit\"at, Mainz, 1976
\end{thebibliography}
\end{document}